\documentclass[journal,twoside,web]{ieeecolor}
\usepackage{generic}
\usepackage{cite}
\usepackage{amsmath,amssymb,amsfonts}
\usepackage{graphicx}
\usepackage{algorithm,algorithmic}
\usepackage{hyperref}
\usepackage{booktabs}
\usepackage{tikz}
\usepackage[caption=false,font=footnotesize]{subfig}
\usepackage[normalem]{ulem}
\hypersetup{hidelinks=true}

\newtheorem{theorem}{Theorem}

\newtheorem{corollary}{Corollary}
\newtheorem{lemma}{Lemma}
\newtheorem{remark}{Remark}
\let\oldremark\remark
\let\endoldremark\endremark

\renewenvironment{remark}
  {\oldremark}
  {\hfill$\square$\endoldremark}

\usetikzlibrary{arrows.meta,calc,backgrounds}

\usepackage{textcomp}
\def\BibTeX{{\rm B\kern-.05em{\sc i\kern-.025em b}\kern-.08em
    T\kern-.1667em\lower.7ex\hbox{E}\kern-.125emX}}
\begin{document}
\title{Policy Iteration over Terminal Cost Functions:\\ From LQR to MPC}
\author{Tren M.J.T. Baltussen, Nathan P. Lawrence, Alexander Katriniok, Ali Mesbah, W.P.M.H. (Maurice) Heemels
\thanks{Tren Baltussen, Alexander Katriniok and Maurice Heemels are with the Department of Mechanical Engineering, Eindhoven University of Technology, The Netherlands. Nathan Lawrence and Ali Mesbah are with the Department of Chemical and Biomolecular Engineering, University of California, Berkeley, USA.
Email: t.m.j.t.baltussen@tue.nl}}

\maketitle

\begin{abstract}
Policy iteration is a classical method for successively improving control policies in optimal control and reinforcement learning. In this paper, we propose a policy iteration framework that refines receding-horizon control (RHC) and model predictive control (MPC) policies online by iteratively updating the terminal cost. We first develop the framework for the linear quadratic regulator (LQR) problem. We show that the Riccati map is a contraction in the Thompson metric on an explicitly constructed invariant order interval. Using this property, we prove convergence of approximate policy iteration for finite-horizon RHC, where the terminal cost is given by an approximate value function. From an initial terminal cost, we then construct a robust invariant order interval within which the RHC policy is refined online. We extend the analysis to constrained MPC and present a constructive algorithm for safely updating the terminal ingredients, which guarantees recursive feasibility and asymptotic stability of the closed loop throughout the iterations. These results connect classical Riccati theory with approximate policy iteration and provide a principled pathway toward learning MPC terminal ingredients. A numerical example illustrates the proposed method.
\end{abstract}

\begin{IEEEkeywords}
Model predictive control (MPC), optimal control, policy iteration, reinforcement learning.
\end{IEEEkeywords}

\section{Introduction}
\label{sec:introduction}

\IEEEPARstart{P}{olicy} iteration (PI) is one of the foundational algorithms
of dynamic programming and model-based reinforcement learning. It
computes an optimal policy by alternating two steps: policy evaluation,
in which the quality of the policy is measured, typically by computing its cost-to-go; and policy improvement, in which a new policy is obtained by optimizing with respect to that cost-to-go. 
Under suitable conditions, this alternation produces a sequence of monotonically improving policies that converges to an optimal one \cite{bertsekas1996neuro}. In its exact form, however, PI is rarely practical. Policy evaluation requires the infinite-horizon cost of a policy over the entire state space, and policy improvement requires an optimization over the space of all admissible policies \cite{LagoudakisParr2003}. 
Any practical PI strategy must therefore rest on two choices: a parameterization of the policies and value functions over which the iteration operates, and approximate
versions of the evaluation and improvement steps that can be carried
out with finite computation and finite data \cite{bertsekas1996neuro}. The properties of the resulting scheme, including whether it converges and whether the intermediate policies are safe to deploy, depend critically on these choices of parameterization and approximation.

In this paper, we argue that receding-horizon control (RHC) and model
predictive control (MPC) provide a particularly well-structured
parameterization for PI. For clarity, we use RHC and MPC to distinguish
between the unconstrained and constrained cases, respectively. In RHC
and MPC, the infinite-horizon optimal control problem is approximated by
a finite-horizon problem with a terminal cost function. The terminal cost function is a standard mechanism to transfer stability and performance properties from an infinite-horizon problem to a finite-horizon formulation \cite{Mayne2000}. When appropriately designed, this terminal cost approximates the value of the infinite-horizon problem beyond the prediction horizon, and can be used to guarantee closed-loop stability and performance \cite{Mayne2000}. After fixing the model, stage cost, horizon,
and constraints, the terminal cost fully determines the resulting control law, and it can therefore be viewed as the parameter of a policy class. The terminal cost thus plays a dual role: it parameterizes the policy and serves as the (approximate) measure of its value.

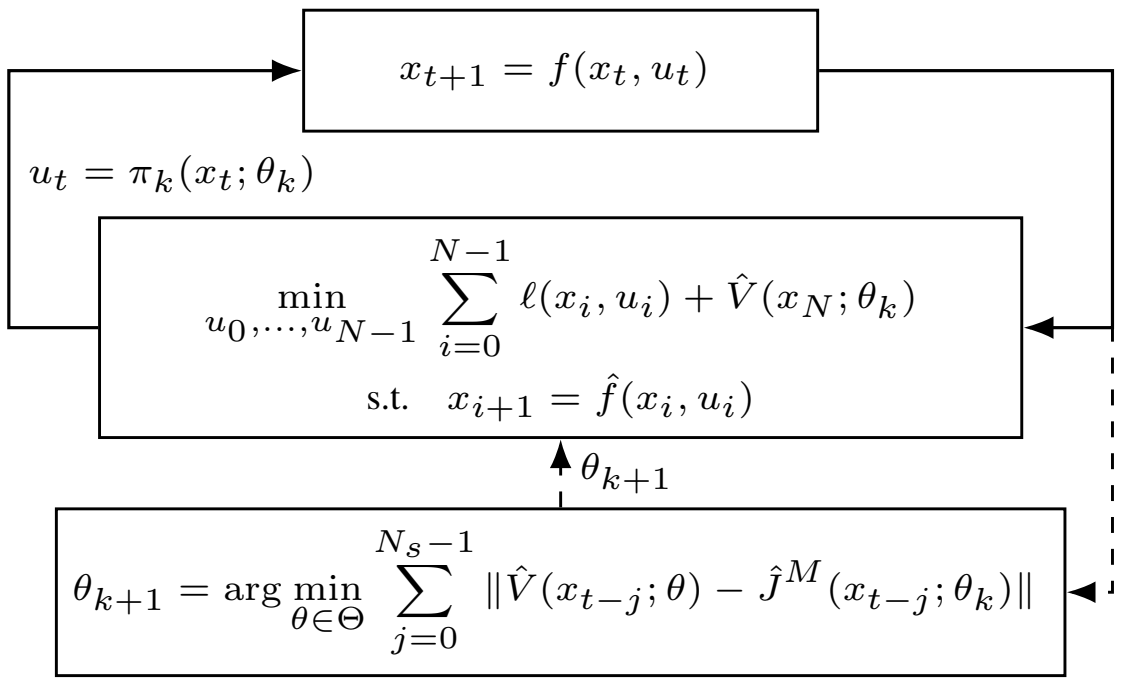
\begin{figure}[t]
\centering
\begin{tikzpicture}[
    >=Latex,
    font=\scriptsize,
    line width=0.55pt,
    block/.style={
        draw,
        fill=white,
        rectangle,
        minimum height=0.8cm,
        align=center,
        inner sep=3pt
    },
    flow/.style={
        ->,
        line width=0.6pt
    }
]


\node[
    block,
    minimum width=3.4cm
] (plant) at (0,1.7)
{
    $\displaystyle x_{t+1}=f(x_t,u_t)$
};

\node[
    block,
    minimum width=6.1cm,
    minimum height=1.45cm
] (mpc) at (0,0)
{
    $\displaystyle
    \min_{u_0,\ldots,u_{N-1}}
    \sum_{i=0}^{N-1}
    \ell(x_i,u_i)+\hat V(x_N;\theta_k)
    $\\[0.7ex]
    $\displaystyle
    \text{s.t.}\quad
    x_{i+1}=\hat f(x_i,u_i)
    $
};

\node[
    block,
    minimum width=4.4cm,
    minimum height=0.8cm
] (update) at (0,-1.75)
{
    $\displaystyle
    \theta_{k+1}
    =
    \arg\min_{\theta\in\Theta}
    \sum_{j=0}^{N_s-1}
    \|
    \hat V(x_{t-j};\theta)
    -
    \hat J^M(x_{t-j};\theta_k)
    \|
    $
};


\coordinate (righttop)  at (3.65,1.7);
\coordinate (rightmid)  at (3.65,0.0);
\coordinate (rightlow)  at (3.65,-1.75);

\coordinate (leftmid)   at (-3.65,0.0);
\coordinate (lefttop)   at (-3.65,1.7);


\draw[flow]
    (plant.east)
    -- (righttop)
    -- node[right] {$x_t$} (rightmid)
    -- (mpc.east);

\draw[flow]
    (mpc.west)
    -- (leftmid)
    -- node[right, pos=0.6] {$u_t=\pi_k(x_t;\theta_k)$} (lefttop)
    -- (plant.west);

\draw[dashed, ->]
    (righttop)
    -- (rightlow)
    -- (update.east);

\draw[dashed, ->]
    (update.north)
    -- node[right] {$\theta_{k+1}$} (mpc.south);

\end{tikzpicture}
\caption{Proposed terminal cost PI architecture.
Closed-loop rollout data from the RHC/MPC policy are used to compute
the next terminal cost parameter $\theta_{k+1}$ using the terminal cost function $\hat V(\cdot;\theta_k)$ as a bootstrap for the finite-horizon evaluation. Here, $t$ denotes plant time (solid lines) and $k$ denotes
the PI update index (dashed lines).}
\label{fig:mpc_architecture}
\end{figure}

We exploit this duality to connect PI to refining RHC and MPC policies
online by iteratively updating the terminal cost, as illustrated in
Fig.~\ref{fig:mpc_architecture}. Given a terminal cost, solving the
finite-horizon optimal control problem constitutes a policy improvement
step, namely an $N$-step lookahead with respect to that terminal cost
rather than the one-step greedy update of classical PI. Deploying the
resulting law in closed loop and accumulating the cost over a finite
window, with the current terminal cost used to bootstrap the remainder,
constitutes an approximate policy evaluation step based on an $M$-step
rollout. The evaluated cost can then be used to update the terminal
cost, which induces a new RHC/MPC law, and the process is repeated
online. Compared with generic function approximators for the policy and value function class, this
parameterization has a distinctive advantage: the policy class inherits
the structural properties of MPC, including explicit constraint
handling and closed-loop stability guarantees derived from the terminal
cost \cite{Banker2025}. This opens the possibility of
performing PI while retaining guarantees for every intermediate policy,
not only for the limit of the PI sequence.

As the finite-horizon planning and the finite-horizon evaluation steps are both approximate, the resulting policy iteration is not covered directly by classical exact PI arguments. Moreover, in stabilization problems a poorly evaluated intermediate policy may fail to be stabilizing and incur infinite cost, while in constrained MPC every deployed policy must additionally preserve feasibility and constraint satisfaction. The central questions are therefore under what conditions terminal cost PI converges, how finite planning $(N)$ and evaluation $(M)$ horizons jointly determine its behavior, and whether closed-loop guarantees can be maintained throughout the iterations.

To investigate the mechanism underlying these questions, we first study
the scheme for the deterministic linear-quadratic regulator (LQR)
problem. We emphasize that the purpose is not to propose a new method
for solving LQR: its optimal solution is available directly from the
discrete-time algebraic Riccati equation (DARE), and classical PI
schemes for LQR are well established \cite{Hewer1971,Kleinman1968}.
Rather, LQR is a setting in which the optimal value class is known,
and terminal cost PI can be analyzed in
closed form. There, the RHC law induced by a terminal cost matrix $P_k$
corresponds to an $N$-fold composition of the \emph{Riccati map}, and the
$M$-step rollout with terminal bootstrap admits a closed-form
representation. This allows us to view truncated policy evaluation as a
perturbation of a contractive Riccati map, and to characterize precisely
how the prediction horizon and rollout length govern convergence.

The LQR analysis is also relevant to the constrained case. When state
and input constraints are enforced, the implicit MPC law is generally
nonlinear and the Riccati structure is no longer globally available.
However, in a neighborhood of the equilibrium where constraints are
inactive, the MPC law coincides with the RHC law of the underlying
unconstrained problem, so the mechanism identified for LQR governs the
local behavior of terminal cost PI for MPC. Building on this connection,
we develop a constructive scheme for updating the terminal cost of an
MPC controller that preserves recursive feasibility, and hence
constraint satisfaction, as well as asymptotic stability of the closed
loop throughout the iterations.

\subsection{Related Work}
\label{sec:related_work}

\subsubsection{Policy Iteration and Terminal Cost Learning}
This work is closest in spirit to the approximation-in-value-space perspective developed by Bertsekas \cite{bertsekas2017dynamic,bertsekas2022lessons}, where multistep lookahead and rollout are combined with an approximation of the cost
beyond the planning horizon. Evaluation by truncated rollout with
terminal bootstrapping is also closely related to optimistic or
modified policy iteration \cite{bertsekas1996neuro}. Here, however, the
value approximation is used specifically as the terminal cost of an
RHC/MPC problem. The policy is therefore induced implicitly by
finite-horizon optimal control rather than updated directly in
feedback-gain space, while the same terminal value is updated from
finite closed-loop rollouts. In this sense, the framework combines RHC
and MPC structure with ideas from reinforcement learning
\cite{Lawrence2025}.
The online refinement of the controller from its own closed-loop data
is also conceptually related to the self-tuning perspective in adaptive
control, in which the controller is updated online to achieve desired
closed-loop properties \cite{Astrom1989}.

Direct policy-optimization methods instead update the feedback gain itself \cite{FazelEtAl2018,HuEtAl2023}.
In uncertain and data-driven LQR settings, exact policy evaluation is
commonly replaced by regression, temporal-difference methods, or
finite-sample estimation; representative examples include least-squares PI and finite-time analyses of approximate PI for LQR \cite{BradtkeEtAl1994,LagoudakisParr2003,KrauthTuRecht2019}.
A complementary model-based approach first identifies the system dynamics from data and then synthesizes a controller while accounting for identification uncertainty
\cite{DeanEtAl2018,DeanTuMatniRecht2019}.

Several learning-based MPC methods use approximate value functions,
learned terminal penalties, or data-dependent terminal ingredients to
improve finite-horizon performance; see, e.g.,
\cite{MorenoMoraEtAl2023,Baltussen2025,Bemporad2024}. Related
constructions appear in iterative learning MPC where terminal costs and sets are computed from previous task executions \cite{RosoliaBorrelli2018}, while \cite{Zuliani2026} studies policy optimization for parametrized MPC policies under model uncertainty using differentiable optimization.

\subsubsection{Riccati Contraction and Order-Theoretic Analysis}
Contraction properties of Riccati maps have been studied in invariant metrics, including the Riemannian analysis of Bougerol \cite{Bougerol1993}, and \emph{Thompson metric} \cite{Thompson1963} results of Lawson and Lim \cite{LawsonLim2007} and Lee and Lim~\cite{LeeLim2008}. Gaubert and Qu~\cite{GaubertQu2014} developed general Thompson-metric contraction formulas for order-preserving maps and applied them to discrete-time Riccati maps, establishing global nonexpansiveness and characterizing when strict contraction holds on the cone of positive definite symmetric matrices. In particular, strict global contraction requires a full-row-rank input matrix, so underactuated systems generally lack a uniform strict contraction on the entire open cone. Krauth, Tu, and Recht~\cite{KrauthTuRecht2019} derive a bounded-set Thompson contraction estimate for the discrete-time Riccati map and use this geometry in their analysis of LQR policy iteration. Related order-theoretic arguments also appear in classical Riccati difference-equation analysis. In particular, Bitmead et al.~\cite{BitmeadEtAl1985} use sub- and supersolution conditions to generate monotone Riccati trajectories and establish stabilizing properties along those trajectories.

Our use of this geometry is motivated by a different requirement: finite-rollout policy evaluation need not preserve the exact comparison relations of classical PI. We construct a domain, an explicit \emph{invariant order interval} on which the nominal Riccati map is contractive. We show that this domain is robust to the perturbation introduced by the finite planning and evaluation horizons, and we use it to establish convergence of the approximate PI scheme.

\subsection{Contributions} 
We study policy iteration on the terminal cost of RHC and MPC, with its theoretical foundation established for the LQR problem. The main contributions are summarized below.

\textbf{1. Terminal cost Policy Iteration Scheme.}
We formulate online terminal cost refinement for RHC and MPC as an
approximate PI scheme in which the terminal cost serves simultaneously
as policy parameter and value approximation. Policy improvement is
performed by $N$-step lookahead through the finite-horizon optimal
control problem, while policy evaluation uses an $M$-step rollout with the approximate value function as a terminal bootstrap.

\textbf{2. Robust Convergence under Finite-Horizon Policy Evaluation.}
For the LQR setting, we analyze the convergence of this scheme by
interpreting truncated policy evaluation as a perturbation of a
contractive Riccati map, and make explicit the roles of $N$ and $M$.
To this end, we construct an invariant domain for the Riccati map that
is robust to the perturbation introduced by the finite-rollout
evaluation and establish an input-to-state stability (ISS)-type contraction bound together with computable conditions for exact convergence of the approximate PI scheme. For
positive-semidefinite state costs, we further show that finite-horizon
observability yields strict contraction of a finite composition of the
Riccati map.

\textbf{3. Extension to Constrained MPC.}
For the constrained setting, we exploit the local LQR structure of MPC to obtain a constructive terminal cost PI scheme. We derive conditions and propose an algorithm for updating the terminal cost and terminal set that ensure recursive feasibility, constraint satisfaction, and asymptotic closed-loop stability throughout the iterations. Numerical results illustrate how the PI viewpoint extends beyond the unconstrained setting while retaining the Riccati contraction as its structural foundation.

\subsection{Notation}
We will use the following notation throughout the paper. $\mathbb{N}_0$ denotes all nonnegative integers, and $\mathbb{R}$ denotes the set of real numbers. We denote by $\mathbb{R}^n$ the $n$-dimensional Euclidean space. 
The set of symmetric matrices is denoted by $\mathbb{S}^n$,  the set of positive definite matrices is denoted by $\mathbb{S}^n_{++}$, and the set of positive semidefinite matrices is denoted by $\mathbb{S}^n_{+}$. 
For a matrix $A\in\mathbb{R}^{n\times n}$, we denote its transpose by $A^\top$. 
For a vector $x\in\mathbb{R}^n$, we denote its Euclidean norm by $\|x\|_2$. 
For a symmetric matrix $P\in\mathbb{S}^n$, we denote its minimum and maximum eigenvalues by $\lambda_{\min}(P)$ and $\lambda_{\max}(P)$, respectively. 
For two symmetric matrices $P,Q\in\mathbb{S}^n$, we write $P\succeq Q$ if $P-Q\succeq 0$, and $P\succ Q$ if $P-Q\succ 0$. The identity matrix of size $n$ is denoted by $I_n$. 

\section{Policy Iteration over Terminal Costs}
\subsection{Finite-Sample Policy Iteration}
\label{sec:sample_based_value_fit}
We aim to build the foundation for a general PI framework for refining RHC/MPC policies online, in closed loop, by refining the terminal cost function using finite rollouts of the currently deployed policy, as illustrated in Fig.~\ref{fig:mpc_architecture}. 
Before we restrict ourselves to the LQR setting, we first outline the general finite-sample PI framework for RHC and MPC.

Let $\hat V(\cdot;\theta)$ denote the terminal value
approximation parameterized by $\theta \in \Theta$. At PI step $k \in \mathbb{N}_0$, the $N$-step planner induces the
RHC policy \(\pi_k(x_t;\theta_k)\) obtained by solving
\begin{equation}
\label{eq:rhc_problem}
\begin{aligned}
\min_{u_0,\ldots,u_{N-1}}
& \sum_{i=0}^{N-1}\ell(x_i,u_i)
+\hat V(x_N;\theta_k) \\
\text{s.t.}\quad & x_{i+1}=\hat f(x_i,u_i) \\
\end{aligned}
\end{equation}
and applying the first input \(u_t = \pi_k(x_t, \theta_k) = u_0^{\star}(x_t;\theta_k),
\) in a receding-horizon fashion. Here, $\hat f$ is a model of the system dynamics, which may be known exactly or identified from data. Constraints are omitted from \eqref{eq:rhc_problem} for brevity. The policy can subsequently be evaluated directly from closed-loop
rollout data. Suppose that $\pi_k$ and $\theta_k$ are held fixed over
an $H$-step rollout.
For each starting state $x_t$, $t=0,1,\ldots,H-M$, define the
$M$-step bootstrapped approximate value sample
\begin{equation}
\label{eq:sample_value_target}
\hat J^M(x_t, \theta_k)
:=
\sum_{\tau=0}^{M-1}
\ell(x_{t+\tau},u_{t+\tau})
+
\hat V(x_{t+M};\theta_k),
\end{equation}
where
\(
u_{t+\tau}
=
\pi_k(x_{t+\tau};\theta_k).
\)
The terminal term \(\hat V(x_{t+M};\theta_k)\) bootstraps the finite rollout with the current value approximation, allowing PI to be initialized and updated using only finite-horizon evaluations.
Importantly, \eqref{eq:sample_value_target} depends only on measured
closed-loop states and inputs, together with the terminal bootstrap; evaluation of the system model is not
required.
An $H$-step rollout therefore provides
\(
N_s:=H-M+1
\)
possibly overlapping value samples
\(
\mathcal D_k
=
\{(x_j,\hat J^M(x_j, \theta_k))\}_{j=0}^{N_s-1}.
\)
The next value approximation can then be obtained, e.g., by minimizing some loss function over the sample set $\mathcal D_k$:
\begin{equation}
\label{eq:general_sample_fit}
\theta_{k+1}
\in
\arg\min_{\theta\in\Theta}
\sum_{j=0}^{N_s-1}
\left\|
\hat V(x_j;\theta)-\hat J^M(x_j, \theta_k)
\right\|.
\end{equation}
Thus, $N$ determines the planning horizon, $M$ the policy-evaluation
horizon, and $H = N_s + M - 1$ the amount of closed-loop data collected before the
value-function update.

We analyze this architecture first for deterministic LQR, whose optimal terminal matrix $P^\star$ solves the DARE, and subsequently extend it to constrained MPC.

\subsection{Deterministic LQR Setting}
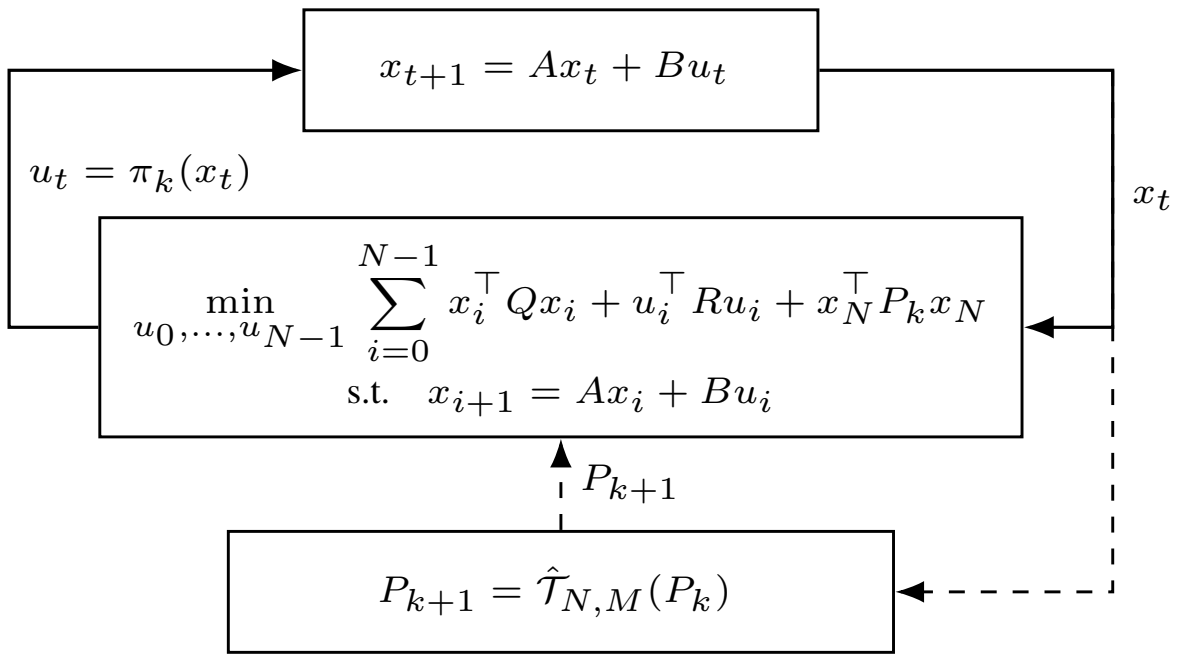
\begin{figure}[t]
\centering
\begin{tikzpicture}[
    >=Latex,
    font=\scriptsize,
    line width=0.55pt,
    block/.style={
        draw,
        fill=white,
        rectangle,
        minimum height=0.8cm,
        align=center,
        inner sep=3pt
    },
    flow/.style={
        ->,
        line width=0.6pt
    }
]


\node[
    block,
    minimum width=3.4cm
] (plant) at (0,1.7)
{
    $\displaystyle x_{t+1}=Ax_t+Bu_t$
};

\node[
    block,
    minimum width=6.1cm,
    minimum height=1.45cm
] (mpc) at (0,0)
{
    $\displaystyle
    \min_{u_0,\ldots,u_{N-1}}
    \sum_{i=0}^{N-1}
    x_i^\top Qx_i+u_i^\top Ru_i+x_N^\top P_k x_N
    $\\[0.7ex]
    $\displaystyle
    \text{s.t.}\quad
    x_{i+1}=Ax_i+Bu_i
    $
};

\node[
    block,
    minimum width=4.4cm,
    minimum height=0.8cm
] (update) at (0,-1.75)
{
    $\displaystyle
    P_{k+1}=\hat{\mathcal T}_{N,M}(P_k)
    $
};


\coordinate (righttop)  at (3.65,1.7);
\coordinate (rightmid)  at (3.65,0.0);
\coordinate (rightlow)  at (3.65,-1.75);

\coordinate (leftmid)   at (-3.65,0.0);
\coordinate (lefttop)   at (-3.65,1.7);


\draw[flow]
    (plant.east)
    -- (righttop)
    -- node[right] {$x_t$} (rightmid)
    -- (mpc.east);

\draw[flow]
    (mpc.west)
    -- (leftmid)
    -- node[right, pos=0.6] {$u_t=\pi_k(x_t)$} (lefttop)
    -- (plant.west);

\draw[dashed, ->]
    (righttop)
    -- (rightlow)
    -- node[below,pos=0.25] {} (update.east); 

\draw[dashed, ->]
    (update.north)
    -- node[right] {$P_{k+1}$} (mpc.south);

\end{tikzpicture}
\caption{Proposed terminal cost PI architecture in the LQR setting.
The RHC policy is parametrized by the terminal cost matrix $P_k$. Here, $t$ denotes plant time (solid lines) and $k$ denotes the PI update index (dashed lines).}
\label{fig:mpc_architecture_lqr}
\end{figure}

As motivated above, we first study the proposed terminal cost PI framework in the deterministic LQR setting, which provides a tractable foundation for analyzing the convergence and closed-loop properties of the proposed scheme. The resulting architecture is illustrated in Fig.~\ref{fig:mpc_architecture_lqr}.
We consider the discrete-time linear system
\[
\label{eq:lin_system}
    x_{t+1} = A x_t + B u_t,
\]
with $x_t \in \mathbb{R}^n$ and $u_t \in \mathbb{R}^m$ at time $t \in \mathbb{N}_0$, with known system matrices $A$ and $B$, and a known stage cost
\[
    \ell(x,u) = x^\top Q x + u^\top R u, \qquad Q \succeq 0, R \succ 0.
\]
We assume that $(A,B)$ is stabilizable and $(A, Q^{1/2})$ is observable. 
We look to minimize the infinite-horizon value function:
\begin{equation}
\label{eq:infinite_horizon_lqr}
J^\infty(x_0) = \min_{u_0,u_1,\ldots} \sum_{t=0}^{\infty} x_t^\top Q x_t + u_t^\top R u_t.
\end{equation}
For the unconstrained LQR setting, the value function is quadratic and the policy is linear.  
As mentioned above, this problem has a well-known solution that is given by the DARE and can be solved iteratively using the Riccati map \cite{bertsekas2017dynamic}.
Here, we assume the system matrices $A$ and $B$, and the cost matrices $Q$ and $R$ are known exactly. Hence, we can define the terminal cost and policy as
\(
\hat V(x;P_k)=x^\top P_kx,
\)
and
\(
\pi_k(x;P_k)=K_N(P_k)x,
\)
respectively.
For every $M$-step rollout window,
\begin{align}
\hat J^M(x_t, P_k)
&=
\sum_{\tau=0}^{M-1}
\ell(x_{t+\tau},K_kx_{t+\tau})
+x_{t+M}^\top P_kx_{t+M},
\end{align}
where $K_k = K_N(P_k)$ is the linear feedback gain induced by the $N$-step RHC law with terminal cost $P_k$.

In general, the value function can be estimated from value samples by fitting the rollout targets to a function in the chosen value-function class, as in \eqref{eq:sample_value_target} and \eqref{eq:general_sample_fit}.
For a quadratic function class, given samples
$\{(x_j,\hat J^M(x_j, P_k))\}_{j=0}^{N_s-1}$, we compute
\begin{equation}
\label{eq:sample_based_lsq}
P_{k+1}
\in
\arg\min_{P\in\mathbb S^n}
\sum_{j=0}^{N_s-1}
\left(
x_j^\top P x_j-\hat J^M(x_j, P_k)
\right)^2,
\end{equation}
which has a unique solution whenever the quadratic regressors are sufficiently exciting,
i.e., the associated regression matrix has rank $d=n(n+1)/2$. 
Hence at least $N_s\ge d$ samples are required. 
If the rank condition is not satisfied, additional trajectories,
initial conditions, or excitation are required.
In the deterministic LQR setting, the terminal cost update \eqref{eq:sample_based_lsq} can be evaluated analytically as detailed in Section~\ref{sec:rhc}.

We aim to address the following questions.\\

\noindent
\fbox{%
\vspace{0.5em}
    \parbox{\dimexpr\columnwidth-2\fboxsep-2\fboxrule\relax}{%
        \textbf{Problem Statement:}
        Can a finite-horizon RHC controller, parameterized by terminal cost \(P_k\),  refine \(P_k\)  online from finite rollouts while preserving closed-loop stability and convergence to the optimal LQR policy? Can this idea be extended to constrained MPC, where the terminal cost and constraint set are updated jointly to maintain recursive feasibility and stability?
    }%
} \\

\section{The Riccati Map}
\label{sec:lqr}
In this setting, finite-horizon optimal control is characterized by repeated application of the Riccati map, whose fixed point corresponds to the infinite-horizon optimal value function. This structure allows us to establish contraction properties that support the subsequent analysis of approximate PI for
RHC and its extension to constrained MPC.

\subsection{The Riccati Map and the Thompson Metric}
Let $J : \mathbb{R}^n \rightarrow \mathbb{R}$ be a value function. The Bellman operator $T$ for the LQR problem \eqref{eq:infinite_horizon_lqr} reads
\[
    (TJ)(x) = \min_u (\ell(x,u) + J(Ax + Bu) ).
\]
In the linear quadratic setting, we can restrict $J$ to be a quadratic function $J_P(x):= x^\top P x$ for $P \in \mathbb{S}_{+}^n$.
Hence,
\[
\begin{aligned}
    & (TJ_P)(x) = \min_u (x^\top Q x + u^\top R u + (x^+)^\top P (x^+) )\\
    & = x^\top (Q + A^\top P A - A^\top P B (R + B^\top P B)^{-1} B^\top P A) x,
\end{aligned}
\]
where $x^+ = Ax + Bu$ and the minimizer is given by {$u^* = -(R + B^\top P B)^{-1} B^\top P A x$.}
So, $T$ maps quadratic functions to quadratic functions, and gives rise to the Riccati map:
\[
    \mathcal{R}(P) := Q + A^\top P A - A^\top P B (R + B^\top P B)^{-1} B^\top P A.
\]
    This formula is well defined on \(\mathbb S_+^n\), since \(R+B^\top PB\succ0\) for every \(P\succeq0\). We restrict to \(\mathbb S_{++}^n\) when introducing the Thompson metric and the inverse representation below.

The policy extraction operator $\mathcal{G}$ returns an optimal policy for one Bellman step through the greedy feedback gain $G(P)$:
\begin{equation}
\label{eq:greedy}
    \mathcal{G}(J_P)(x) = G(P) x, \quad G(P) = -(R + B^\top P B)^{-1} B^\top P A.
\end{equation}
From the DARE, we obtain the fixed-point condition
\[
    J_{P^\star} = T(J_{P^\star}) \iff P^\star = \mathcal{R}(P^\star).
\]
This fixed point, sometimes referred to as the steady-state solution, can be obtained by iterating the Riccati map backward starting from any terminal condition $P_0 \in \mathbb{S}_{++}^n$ \cite{bertsekas2017dynamic}.
By establishing contraction properties of the Riccati map, we can analyze the convergence of the terminal cost PI scheme through Banach's fixed-point theorem \cite{Bullo2026}.

Before we proceed analyzing $\mathcal{R}$, we first
recall the order and metric structures used throughout the paper.
The vector space $\mathbb S^n$ of symmetric matrices is partially
ordered by the Loewner order,
\(
X\preceq Y
\iff 
Y-X\in\mathbb S_+^n.
\)
Given $P_-,P_+\in\mathbb S_{++}^n$ with $P_-\preceq P_+$, the
corresponding \emph{order interval} is
\(
[P_-,P_+]
:=
\{P\in\mathbb S_{++}^n:
P_-\preceq P\preceq P_+\}.
\)
An order interval $\mathcal P$ is said to be invariant under a map
$\mathcal R$ if
\(
\mathcal R(\mathcal P)\subseteq\mathcal P.
\)
On the interior $\mathbb S_{++}^n$ of the positive-semidefinite cone,
we consider the \emph{Thompson metric}~\cite{Thompson1963},
\begin{equation}
\label{eq:Thompson_metric}
d_T(X,Y)
=
\log\max\{
\lambda_{\max}(X^{-1}Y),
\lambda_{\max}(Y^{-1}X)
\},
\end{equation}
for $X,Y\in\mathbb S_{++}^n$.
The metric space $(\mathbb S_{++}^n,d_T)$ is complete
\cite{Bullo2026} and the Thompson metric has been used in the analysis
of Riccati recursions
\cite{LawsonLim2007,GaubertQu2014}. It is invariant under invertible
congruence transformations and is naturally compatible with the cone
order, making it well suited to the order-preserving structure of the
Riccati map.

As the Riccati map is not globally contractive \cite{GaubertQu2014}, we study its
contraction properties on an invariant order interval
$\mathcal P\subset\mathbb S_{++}^n$. Specifically, we seek
$\gamma\in[0,1)$ such that
\[
d_T(\mathcal R(X),\mathcal R(Y))
\le
\gamma d_T(X,Y),
\qquad
X,Y\in\mathcal P.
\]

\subsection{Constructing an Invariant Order Interval}
To apply Banach's fixed-point theorem, we construct an order interval
$\mathcal P=[P_-,P_+]\subset\mathbb S_{++}^n$ such that
$\mathcal R(\mathcal P)\subseteq\mathcal P$ and
$(\mathcal P,d_T)$ is complete. We first establish monotonicity of the
Riccati map in the Loewner order and then construct suitable lower and
upper endpoints.
For strict contractivity of $\mathcal{R}$, we require $Q \succ 0$, i.e., we need $P$ to be uniformly bounded away from the boundary \cite{Bullo2026}. However, by considering a finite composition of the Riccati map, we can generalize the results to $Q \succeq 0$.\\

\begin{lemma}[Self-Map]
    \label{lem:riccati_self_map}
    If $Q \succ 0, R \succ 0$, then the Riccati map is a self-map, i.e., $\mathcal{R} : \mathbb{S}^n_{++} \rightarrow \mathbb{S}^n_{++}$. 
\end{lemma}

\begin{proof}
    The Riccati map can be expressed as
    \[
    \mathcal{R}(P) = \min_K (Q + K^\top R K + (A+BK)^\top P (A+BK)).
    \]
    Since \[(A + BK)^\top P (A + BK) \succeq 0, \quad K^\top R K \succeq 0  \quad \forall K,\]
    we have \( \mathcal{R}(P) \succeq Q \succ 0. \)
    Hence, $\mathcal{R}(P) \in \mathbb{S}^n_{++}$.
\end{proof}
\,

\begin{lemma}[Order Preservation of the Riccati map]
\label{lem:riccati_mono}
    If $P \preceq P'$ for $P,P' \in \mathbb{S}_{++}^n$, then
    \(\mathcal{R}(P) \preceq \mathcal{R}(P').\)
\end{lemma}

\begin{proof}
    For any fixed feedback gain $K \in \mathbb{R}^{m \times n}$ we define \[ \Phi(P,K) := Q + K^\top R K + (A+BK)^\top P (A + BK).\]
    For the LQR, we have \(\mathcal{R}(P) = \min_K \Phi(P,K)\) with the greedy feedback gain
    \[G(P) = \arg \min_K \Phi(P,K) = -(R + B^\top P B)^{-1} B^\top P A.\]

    Let \(K':=G(P')\). Since \(P\preceq P'\), \[ \Phi(P,K') \preceq \Phi(P',K'). \] Moreover, by optimality of \(G(P)\), \[ \mathcal R(P) = \min_K \Phi(P,K) \preceq \Phi(P,K'). \] Hence \[ \mathcal R(P) \preceq \Phi(P,K') \preceq \Phi(P',K') = \mathcal R(P'), \] which proves that \(\mathcal R\) is order preserving.
\end{proof}
\,

Let us define the Riccati \emph{subsolution} and \emph{supersolution} conditions, respectively, as:
\begin{align}
    \mathcal{R}(P_{-}) &\succeq P_{-}, \\
    \mathcal{R}(P_{+}) &\preceq P_{+}.
\end{align}

In order to construct an upper bound on the interval $\mathcal{P}$, consider any stabilizing gain
$\bar K$ such that
\(
    A_{\bar K}:=A+B\bar K
\)
is Schur. Define
\(
    S_K:=Q+K^\top RK\succ0.
\)
The corresponding policy-evaluation matrix $P_{\bar K}\succ0$ is
\[
    P_{\bar K}
    =
    \sum_{t=0}^{\infty}
    (A_{\bar K}^t)^\top
    S_{\bar K}
    A_{\bar K}^t,
\]
which is the unique solution of the Lyapunov equation
\(
    P_{\bar K}
    =
    S_{\bar K}
    +
    A_{\bar K}^\top P_{\bar K}A_{\bar K}.
\)
Since the Riccati map minimizes the quadratic Bellman update over the
feedback gain,
\begin{equation}
    \label{eq:supersolution}
    \mathcal R(P_{\bar K})
    \preceq
    S_{\bar K}
    +
    A_{\bar K}^\top P_{\bar K}A_{\bar K}
    =
    P_{\bar K},
\end{equation}
so $P_{\bar K}$ is a supersolution.
Moreover, the supersolution can be enlarged arbitrarily. Let
$Z_{\bar K}\succ0$ denote the unique solution of
\(
    Z_{\bar K}
    -
    A_{\bar K}^\top Z_{\bar K}A_{\bar K}
    =
    I,
\)
and define
\[
P_{\max}(\beta)
:=
P_{\bar K}+\beta Z_{\bar K},
\qquad \beta\ge0.
\]
Then
\[
\begin{aligned}
\mathcal R(P_{\max}(\beta))
&\preceq
S_{\bar K}
+A_{\bar K}^\top P_{\max}(\beta)A_{\bar K}\\
&=
P_{\max}(\beta)-\beta I
\preceq P_{\max}(\beta).
\end{aligned}
\]
Hence $P_{\max}(\beta)$ is a supersolution for every
$\beta\ge0$. Because $Z_{\bar K}\succ0$, these
supersolutions can be made arbitrarily large in the
Loewner order by increasing $\beta$.

We therefore define the order interval
\[
    \mathcal P
    :=
    \{P\in\mathbb S_{++}^n:
      Q\preceq P\preceq P_{\max}\},
\]
where $P_{\max}$ denotes any such supersolution.

\begin{lemma}[Invariant Order Interval]
\label{lem:riccati_inv}
Let \(P_{\max}\succ0\) satisfy
\[
    \mathcal R(P_{\max})\preceq P_{\max}.
\]
Then the order interval \(\mathcal{P}=[Q,P_{\max}]\) is invariant under
the Riccati map \(\mathcal{R}\).
\end{lemma}

\begin{proof}
For \(P\in\mathcal P\), Lemma~\ref{lem:riccati_mono} and the
supersolution condition give
\[
    \mathcal R(P)\preceq\mathcal R(P_{\max})\preceq P_{\max}.
\]
Moreover, \(\mathcal R(P)\succeq Q\) for every \(P\succeq0\).
Therefore \(Q\preceq\mathcal R(P)\preceq P_{\max}\), proving
\(\mathcal R(P)\in\mathcal P\).
\end{proof}

\begin{lemma}[Completeness in the Thompson metric]
\label{lem:riccati_complete}
    The order interval $\mathcal{P}$ is complete under the Thompson metric $d_T$.
\end{lemma}
\begin{proof}
The cone \((\mathbb S^n_{++},d_T)\) is complete. Moreover, on sets bounded
away from the boundary of \(\mathbb S^n_{++}\), the Thompson topology agrees
with the usual matrix topology. Since \(Q\succ0\), the interval
\(\mathcal P=[Q,P_{\max}]\) is bounded away from the boundary, and it is
closed in the usual matrix topology. Hence it is closed in the Thompson
topology. Therefore \((\mathcal P,d_T)\) is a closed subset of a complete
metric space, and is complete.
\end{proof}

\subsection{Semi-global Contraction of the Riccati Map}
We now establish a uniform contraction estimate on the invariant order
interval constructed above. The proof builds on the translation-contraction
property of invariant metrics developed by Lee and Lim~\cite{LeeLim2008}.
By exploiting the invariance of $\mathcal P$, we obtain a
uniform bound on the nonconstant part of the Riccati map relative to $Q$,
yielding a contraction certificate tailored to the subsequent
perturbation analysis.

\begin{theorem}
\label{thm:contraction}
Let \(Q\succ0\), \(R\succ0\), and let \(P_{\max}\succ0\) be a supersolution
of the Riccati map, i.e.,
\(
    \mathcal R(P_{\max})\preceq P_{\max}.
\)
Define \(\mathcal P=[Q,P_{\max}]\),
\(
    \Delta_{\mathcal P}:=d_T(Q,P_{\max}),
\)
and
\[
    \Lambda_{\mathcal P}
    :=
    \lambda_{\max}
    \left(Q^{-1/2}P_{\max}Q^{-1/2}\right)
    =
    e^{\Delta_{\mathcal P}}.
\]
Then $(\mathcal P,d_T)$ is a complete metric space,
\[
    \mathcal R(\mathcal P)\subseteq\mathcal P,
\]
and $\mathcal R:\mathcal P\to\mathcal P$ is a contraction satisfying
\[
    d_T(\mathcal R(P),\mathcal R(P'))
    \le
    \bar\gamma_{\mathcal P}\,
    d_T(P,P'),
    \qquad P,P'\in\mathcal P,
\]
where
\[
    \boxed{
    \bar\gamma_{\mathcal P}
    :=
    1-\Lambda_{\mathcal P}^{-1}
    =
    1-e^{-\Delta_{\mathcal P}}
    \in[0,1).
    }
\]

Consequently, \(\mathcal R\) admits a unique fixed point
\(P^\star\in\mathcal P\), and for every \(P_0\in\mathcal P\) the Picard
iteration \(P_{k+1}=\mathcal R(P_k)\), with $k\in\mathbb N$, satisfies
\[
    d_T(P_k,P^\star)
    \le
    \bar\gamma_{\mathcal P}^{\,k}
    d_T(P_0,P^\star),
    \quad \text{for all } k \in \mathbb N.
    \vspace{0.5em}
\]
\end{theorem}

\begin{proof}
By Lemmas~\ref{lem:riccati_self_map}--\ref{lem:riccati_complete},
$\mathcal R$ is order preserving,
$\mathcal R(\mathcal P)\subseteq\mathcal P$, and
$(\mathcal P,d_T)$ is complete. It remains to establish the contraction
estimate.

Using the Woodbury identity, write
\[
    \mathcal R(P)=Q+Z(P),
    \;
    Z(P):=
    A^\top
    (P^{-1}+BR^{-1}B^\top)^{-1}
    A
    \succeq0.
\]
Let $P_1,P_2\in\mathcal P$ and set
$\delta:=d_T(P_1,P_2)$ and $\zeta:=e^\delta$.
By the definition of the Thompson metric,
\[
    P_1\preceq\zeta P_2,
    \qquad
    P_2\preceq\zeta P_1.
\]
From $P_1\preceq\zeta P_2$ and the order-reversing property of matrix
inversion, $P_1^{-1}\succeq\zeta^{-1}P_2^{-1}$.
Since $BR^{-1}B^\top\succeq0$ and $\zeta\ge1$,
\[
\begin{aligned}
    P_1^{-1}+BR^{-1}B^\top
    &\succeq
    \zeta^{-1}P_2^{-1}+BR^{-1}B^\top\\
    &\succeq
    \zeta^{-1}
    (P_2^{-1}+BR^{-1}B^\top).
\end{aligned}
\]
Inverting once more and applying congruence by $A$ gives
\(
    Z(P_1)\preceq\zeta Z(P_2).
\)
Interchanging $P_1$ and $P_2$ yields
\(
    Z(P_2)\preceq\zeta Z(P_1).
\)

Next, invariance of $\mathcal P$ gives
\(
    Q+Z(P)=\mathcal R(P)\preceq P_{\max},
\)
for every $P\in\mathcal P$.
By definition of $\Lambda_{\mathcal P}$,
\(
    P_{\max}\preceq\Lambda_{\mathcal P}Q,
\)
and hence, for all \(P\in\mathcal P\),
\begin{equation}
\label{eq:subhomogen}
    0\preceq Z(P)
    \preceq(\Lambda_{\mathcal P}-1)Q.
\end{equation}

To satisfy the positive-definiteness requirement of the translation
inequality in \cite[Prop.~4.1]{LeeLim2008}, define, for
$\varepsilon>0$,
\[
    Z_i^\varepsilon:=Z(P_i)+\varepsilon Q,
    \qquad i=1,2.
\]
Then $Z_i^\varepsilon\succ0$. Moreover, the preceding order relations imply
$Z_1^\varepsilon\preceq\zeta Z_2^\varepsilon$ and
$Z_2^\varepsilon\preceq\zeta Z_1^\varepsilon$, so that
$d_T(Z_1^\varepsilon,Z_2^\varepsilon)\le\delta$.
Furthermore,
\[
    Z_i^\varepsilon
    \preceq
    (\Lambda_{\mathcal P}-1+\varepsilon)Q.
\]

By congruence invariance of the Thompson metric, define
\[
    \widetilde Z_i^\varepsilon
    :=
    Q^{-1/2}Z_i^\varepsilon Q^{-1/2}.
\]
Then
\[
    0\prec\widetilde Z_i^\varepsilon
    \preceq
    (\Lambda_{\mathcal P}-1+\varepsilon)I.
\]
Applying \cite[Prop.~4.1]{LeeLim2008} to the translation by $I$ gives
\[
\begin{aligned}
    d_T(Q+Z_1^\varepsilon,Q+Z_2^\varepsilon)
    &=
    d_T(I+\widetilde Z_1^\varepsilon,
        I+\widetilde Z_2^\varepsilon)\\
    &\le
    \frac{\Lambda_{\mathcal P}-1+\varepsilon}
         {\Lambda_{\mathcal P}+\varepsilon}
    d_T(\widetilde Z_1^\varepsilon,
        \widetilde Z_2^\varepsilon)\\
    &\le
    \frac{\Lambda_{\mathcal P}-1+\varepsilon}
         {\Lambda_{\mathcal P}+\varepsilon}\,
    \delta .
\end{aligned}
\]
Letting $\varepsilon\downarrow0$ and using continuity of $d_T$ on
$\mathbb S_{++}^n$ yields
\[
    d_T(\mathcal R(P_1),\mathcal R(P_2))
    \le
    \frac{\Lambda_{\mathcal P}-1}
         {\Lambda_{\mathcal P}}
    d_T(P_1,P_2).
\]
Therefore,
\[
    \bar\gamma_{\mathcal P}
    =
    1-\Lambda_{\mathcal P}^{-1}.
\]
Since $Q\preceq P_{\max}$,
\[
    \Delta_{\mathcal P}
    =
    d_T(Q,P_{\max})
    =
    \log\Lambda_{\mathcal P},
\]
and hence
\[
    \bar\gamma_{\mathcal P}
    =
    1-e^{-\Delta_{\mathcal P}}
    \in[0,1).
\]

Finally, Banach's fixed-point theorem implies that \(\mathcal R\) has a
unique fixed point \(P^\star\in\mathcal P\), and for every
\(P_0\in\mathcal P\),
\[
    d_T(P_k,P^\star)
    \le
    \bar\gamma_{\mathcal P}^{\,k}
    d_T(P_0,P^\star),
    \qquad
    P_{k+1}=\mathcal R(P_k).
    \vspace{-1em}
\]
\end{proof}

The result is semi-global in the sense that the certified invariant
domain can be enlarged arbitrarily through the choice of the
supersolution $P_{\max}$. However, as the interval grows, the contraction factor deteriorates, with \(\bar\gamma_{\mathcal P}\to1\) as \(d_T(Q,P_{\max})\to\infty\).

Related Thompson-metric bounds have been used in the LQR policy-iteration analysis of Krauth, Tu, and Recht~\cite{KrauthTuRecht2019}, where the nonconstant Riccati term $Z$ is controlled through the relaxation
\[
    A^\top(P^{-1}+BR^{-1}B^\top)^{-1}A
    \preceq A^\top P A.
\]
Here, instead, invariance of $\mathcal P$ gives directly
\[
    \mathcal R(P)-Q
    \preceq P_{\max}-Q
    \preceq(\Lambda_{\mathcal P}-1)Q,
\]
which yields
\(
    \bar\gamma_{\mathcal P}=1-\Lambda_{\mathcal P}^{-1}.
\)
Consequently, once the invariant interval has been constructed, $\bar\gamma_{\mathcal P}$ depends only on its relative Thompson geometry. The two contraction estimates are generally incomparable. Our main advantage is therefore not a uniformly tighter contraction factor, but that the same invariant order interval subsequently provides the robustness domain for the truncated policy-evaluation update introduced in Section~\ref{sec:rhc}. In particular, the finite-rollout error can be compared directly against the available margin to the boundary of $\mathcal P$, allowing the approximate update to be analyzed as a perturbation of a contractive self-map on this common invariant domain.

\begin{remark}
The regional nature of the estimate also permits a progressively tighter
contraction certificate along supersolutions of the Riccati equation. In particular, if
\(P_0=P_{\max}\) and
\(
    P_{k+1}=\mathcal R(P_k),
\)
then
\(
    P^\star\preceq P_{k+1}\preceq P_k,
\)
so each interval
\(
    \mathcal P_k=[Q,P_k]
\)
is positively invariant and nested:
\(
    \mathcal P_{k+1}\subseteq\mathcal P_k.
\)
Applying the same contraction estimate on \(\mathcal P_k\) yields the
iteration-dependent certificate
\[
    \bar\gamma_k
    =
    1-\exp\!\left(-d_T(Q,P_k)\right).
\]
Since \(P_k\downarrow P^\star\), the interval diameter decreases to
\(d_T(Q,P^\star)\), and hence the certified contraction factor approaches
the corresponding factor on the limiting interval \([Q,P^\star]\).
Here, $P_k\downarrow P^\star$ denotes monotone convergence in the
Loewner order, i.e.,
\(
P^\star\preceq P_{k+1}\preceq P_k,
\) such that
\(
P_k\to P^\star.
\)
Conversely, initialization from a Riccati subsolution
$P_0\preceq\mathcal R(P_0)$ yields monotone convergence from below,
denoted by $P_k\uparrow P^\star$, provided the iterates remain in the
invariant interval.
\end{remark}

\subsection{Positive Semidefinite Stage Cost}
The previous results can be extended to the case of a positive semidefinite state cost matrix \(Q\succeq0\) by considering a finite composition of the Riccati map. This allows us to establish a contraction property on an invariant order interval for the finite Riccati recursion when \((A,Q^{1/2})\) is observable.

\begin{corollary}[Contraction for semidefinite state costs]
\label{cor:finite_step_contraction}
Assume \(Q\succeq0\), \(R\succ0\), and suppose that
\((A,Q^{1/2})\) is observable. Let $P_{\max}\succ0$ be a supersolution of the Riccati map, i.e.,
\(
    \mathcal R(P_{\max})\preceq P_{\max}.
\)
Then there exists \(\nu\le n\) such that
\(
    S_\nu:=\mathcal R^\nu(0)\succ0.
\)
For any such \(\nu\), define
\(
    \mathcal P_\nu:=[S_\nu,P_{\max}].
\)
Then \((\mathcal P_\nu,d_T)\) is complete,
\(
    \mathcal R^\nu(\mathcal P_\nu)\subseteq\mathcal P_\nu,
\)
and \(\mathcal R^\nu\) is a strict Thompson-metric contraction on
\(\mathcal P_\nu\).

More precisely, let
\(
    \Delta_\nu:=d_T(S_\nu,P_{\max}),
\)
then for all $P,P'\in\mathcal P_\nu$,
\(
    d_T\!\left(
        \mathcal R^\nu(P),\mathcal R^\nu(P')
    \right)
    \le
    \bar\gamma_\nu\,d_T(P,P'),
\)
where
\[
    \boxed{
    \bar\gamma_\nu
    :=
    1-e^{-\Delta_\nu}
    \in[0,1).
    }
\]

Consequently, \(\mathcal R^\nu\) has a unique fixed point
\(P^\star\in\mathcal P_\nu\), which is also a fixed point of
\(\mathcal R\). Moreover, for every \(P_0\in\mathcal P_\nu\),
\[
    (\mathcal R^\nu)^k(P_0)\to P^\star,
\]
with
\(
    d_T\!\left((\mathcal R^\nu)^k(P_0),P^\star\right)
    \le
    \bar\gamma_\nu^{\,k}
    d_T(P_0,P^\star), \text{for all } k \in \mathbb N.
\)
\end{corollary}

\begin{proof}
Observability implies that, for some \(\nu\le n\), the finite-horizon
value matrix \(S_\nu=\mathcal R^\nu(0)\) is positive definite.
Set $\Phi:=\mathcal R^\nu$. Since $\mathcal R$ is monotone and
$P_{\max}$ is a supersolution,
\[
    \Phi(S_\nu)\succeq S_\nu,
    \qquad
    \Phi(P_{\max})\preceq P_{\max},
\]
so $\mathcal P_\nu=[S_\nu,P_{\max}]$ is invariant under $\Phi$ and
complete in $d_T$.

It remains only to recover the strict subhomogeneity estimate from \eqref{eq:subhomogen} in Theorem~\ref{thm:contraction}.  Evaluating the problem
with terminal matrix $\zeta P$ along the same input sequence gives
\[
\begin{aligned}
x^\top\Phi(\zeta P)x
&\le
\sum_{i=0}^{\nu-1}\ell(x_i,u_i)
+\zeta x_\nu^\top P x_\nu \\
&=
\zeta x^\top\Phi(P)x
-(\zeta-1)
\sum_{i=0}^{\nu-1}\ell(x_i,u_i).
\end{aligned}
\]
Since the accumulated stage cost along any input sequence is bounded
below by the zero-terminal $\nu$-step value,
\[
\sum_{i=0}^{\nu-1}\ell(x_i,u_i)
\ge
x^\top S_\nu x,
\]
and therefore
\(
    \Phi(\zeta P)
    \preceq
    \zeta\Phi(P)-(\zeta-1)S_\nu.
\)
Repeating the translation argument of Theorem~1 for
$\Phi=\mathcal R^\nu$ with $S_\nu$ in place of $Q$ therefore yields
\(
    d_T(\mathcal R^\nu(P),\mathcal R^\nu(P'))
    \le
    \left(1-e^{-\Delta_\nu}\right)d_T(P,P'),
\)
which proves the stated contraction factor and convergence.

Finally,
\(
    \mathcal R(S_\nu)\succeq S_\nu,
\)
and
\(
    \mathcal R(P_{\max})\preceq P_{\max},
\)
so $\mathcal R(\mathcal P_\nu)\subseteq\mathcal P_\nu$.
If $P^\star$ is the unique fixed point of $\mathcal R^\nu$, then
$\mathcal R(P^\star)$ is also a fixed point of $\mathcal R^\nu$ in
$\mathcal P_\nu$. Hence, by uniqueness,
\[
    \mathcal R(P^\star)=P^\star.
    \vspace{-1em}
\]
\end{proof}
\,

\begin{remark}[Detectability of the pair \((A,Q^{1/2})\)]
Under detectability alone, stable unobservable directions may prevent $\mathcal{R}^\nu(0)$ from becoming positive definite for any finite $\nu$. Extending the contraction result to this case therefore requires a metric defined on the corresponding cone face or quotient space, rather than the standard Thompson metric on $\mathbb{S}^n_{++}$.
\end{remark}

\section{Receding Horizon Control}
\label{sec:rhc}
Receding-horizon control (RHC) approximates the infinite-horizon problem by repeatedly solving a finite-horizon optimal control problem. 
In this section, we use the contraction properties of the Riccati map on the invariant order interval, constructed in the previous section, to study the convergence of a truncated policy evaluation of the RHC policy.

\subsection{The Finite-Horizon LQR Problem}
We consider the RHC policy constituted by the $N$-step look-ahead minimization with some terminal cost matrix $P$,
\[
    J^N_{P}(x_0) = \min_{u_0, \ldots, u_{N-1}} \sum_{i=0}^{N-1} x_i^\top Q x_i + u_i^\top R u_i + x_N^\top P x_N,
\]
and applying the first action in a receding-horizon manner.
Then, the $N$-step look-ahead minimization of the RHC policy corresponds to iterating $\mathcal{R}$, i.e., 
\(
    T^N(J_{P}) = J_{\mathcal{R}^N(P)}.
\)
Now, we can rewrite the RHC value function as
\[
    J^N_{P}(x) = x^\top P_N(P) x = x^\top \mathcal{R}^N(P) x,
\]
which is the \textit{open-loop} optimum of the RHC, i.e., the RHC control law minimizes an $N$-fold composition of the exact Riccati map.
Recall that the induced RHC policy is given by
\[
    \pi^N_{P}(x) = K_N(P)x = G(\mathcal{R}^{N-1}(P))x,
\]
where $G$ is the greedy feedback \eqref{eq:greedy}.
Here, the first action is greedy with respect to the \(N-1\)-step look-ahead value \cite{bertsekas2017dynamic}. 
Note that both operators have the same fixed point $P^\star$ from the DARE.
Hence, if the terminal cost matrix satisfies $P=P^\star$, then the RHC
law coincides with the infinite-horizon optimal LQR policy, a well-known
result in the RHC literature \cite{Mayne2000}.

\subsection{Truncated Policy Evaluation}
We extend the contraction result of the exact Riccati map by relying on a finite-horizon look-ahead minimization, and a finite-horizon policy evaluation, i.e., a \textit{truncated policy evaluation}.
Given a terminal matrix $P$, we can evaluate the first RHC gain $K_N(P)$, and then evaluate the resulting policy for $M$ time steps with a bootstrapped terminal value matrix $P$.
Then, for the truncated policy evaluation, we evaluate the fixed gain $K$ for $M$ time steps with a bootstrapped terminal matrix $P$:
\begin{equation}
    \hat{P}(K, P, M) = \sum_{t = 0}^{M-1} (A_K^t)^\top S_K A_K^t + (A_K^M)^\top P A_K^M,
\end{equation}
where $A_K = A+BK$ and $S_K = Q + K^\top R K$.
Now, we define the truncated RHC policy iteration operator $\hat{\mathcal{T}}_{N,M}$ as
\begin{equation}
\label{eq:truncated_pi_update}
    \hat{\mathcal{T}}_{N,M}(P) := \hat{P}(K_N(P), P, M),
\end{equation}
and study its convergence properties as a perturbation of the composed Riccati map $\mathcal{R}^{N}$, cf. Fig.~\ref{fig:mpc_architecture_lqr}.\\

\noindent
\fbox{%
\vspace{0.5em}
    \parbox{\dimexpr\columnwidth-2\fboxsep-2\fboxrule\relax}{%
        \textbf{Approximate Policy Iteration Algorithm:}

        \noindent
        Given $N \geq 1, M \geq 1$, and $P_0$. For $k=0,1,2,\ldots$:

\noindent
\textbf{1. Policy improvement} via $N$-step lookahead minimization.
    \[
        \label{eq:rhc_policy_improvement}
        \pi_k = \mathcal{G}(T^{N-1} \hat{J}_k) \iff K_k = K_N(P_k) = G(\mathcal{R}^{N-1}(P_k)).
    \]

\noindent
\textbf{2. Policy evaluation} via truncated application of $T_{\pi_k}$ with approximate terminal value function $\hat{J_k}$. 
    \[
        \label{eq:rhc_policy_evaluation}
        \hat{J}_{k+1}^M = (T_{\pi_k})^M \hat{J}_k \iff P_{k+1} = \hat{\mathcal{T}}_{N,M}(P_k).
    \]
    }%
}\\

\subsection{ISS-Type Contraction of Truncated Policy Iteration}
\label{sec:truncated_pi_contraction}

We next study convergence of
\[
    P_{k+1}
    =
    \hat{\mathcal T}_{N,M}(P_k),
    \qquad
    P_k\in\mathcal P_\nu=[S_\nu,P_{\max}],
\]
by viewing the truncated PI map as a finite-horizon perturbation of the
contractive Riccati planner $\mathcal R^N$. 
The truncated-PI operator shares the Riccati fixed point
$P^\star$, since $\hat{\mathcal T}_{N,M}(P^\star)=P^\star$
for all $N,M\ge1$. We first derive an
a-priori certificate ensuring that $\mathcal P_\nu$ is invariant under
the truncated PI update. We then establish an ISS-type contraction bound \cite{sontag1989} and combine it with the local structure of the truncated PI map to obtain computable conditions for exact convergence to the unique fixed point $P^\star$ of the Riccati map.

For the robust-invariance certificate below, take
$P_{\max}$ to be a strict Riccati supersolution and define
$\eta:=\lambda_{\min}(P_{\max}-\mathcal R(P_{\max}))>0$.
Then
\[
    \mathcal R(P_{\max})\preceq P_{\max}-\eta I.
\]
Since the Riccati iterates from $P_{\max}$ are decreasing, for every
$P\in\mathcal P_\nu$ and $M\ge1$,
\begin{equation}
\label{eq:strict_upper_slack}
    \mathcal R^M(P)
    \preceq
    \mathcal R^M(P_{\max})
    \preceq
    \mathcal R(P_{\max})
    \preceq
    P_{\max}-\eta I.
\end{equation}

\begin{lemma}[Finite-horizon perturbation and invariance]
\label{lem:robust_pi_invariance}
For $P\in\mathcal P_\nu$, let
\[
    K:=K_N(P)=G(\mathcal R^{N-1}(P)),
    \qquad
    A_K:=A+BK,
\]
and, for $i=0,\ldots,M-1$, define
\[
    \Pi_i(P)
    :=
    \bigl(K-G(\mathcal R^i(P))\bigr)^\top
    \bigl(R+B^\top\mathcal R^i(P)B\bigr)
    \bigl(K-G(\mathcal R^i(P))\bigr).
\]
Then
\begin{equation}
\label{eq:finite_pi_perturbation}
    \hat{\mathcal T}_{N,M}(P)-\mathcal R^N(P)
    =
    \mathcal R^M(P)-\mathcal R^N(P)
    +
    C_{N,M}(P),
\end{equation}
where
\begin{equation}
\label{eq:C_NM}
    C_{N,M}(P)
    :=
    \sum_{i=0}^{M-1}
    (A_K^{M-1-i})^\top
    \Pi_i(P)
    A_K^{M-1-i}
    \succeq0.
\end{equation}

Define
\[
    D_M(P):=P_{\max}-\mathcal R^M(P).
\]
Then $D_M(P)\succeq\eta I$ on $\mathcal P_\nu$. Moreover,
\begin{equation}
\label{eq:rho_NM}
    \rho_{N,M}
    :=
    \max_{P\in\mathcal P_\nu}
    \lambda_{\max}\!\left(
        D_M(P)^{-1/2}
        C_{N,M}(P)
        D_M(P)^{-1/2}
    \right)
\end{equation}
is finite and attained, and
\begin{equation}
\label{eq:rho_invariance}
\boxed{
    \hat{\mathcal T}_{N,M}(\mathcal P_\nu)
    \subseteq
    \mathcal P_\nu
    \quad\Longleftrightarrow\quad
    \rho_{N,M}\le1.
}
\vspace{1em}
\end{equation}
\end{lemma}

\begin{proof}
For fixed $K=K_N(P)$, completion of squares gives
\[
    Q+K^\top RK
    +A_K^\top\mathcal R^i(P)A_K
    =
    \mathcal R^{i+1}(P)+\Pi_i(P).
\]
With the convention $\hat P(K,P,0)=P$, define
\[
    E_m:=\hat P(K,P,m)-\mathcal R^m(P).
\]
Then $E_0=0$ and
\(
    E_{m+1}
    =
    A_K^\top E_mA_K+\Pi_m(P).
\)
Iteration yields
\[
    E_M
    =
    \sum_{i=0}^{M-1}
    (A_K^{M-1-i})^\top
    \Pi_i(P)
    A_K^{M-1-i}
    =
    C_{N,M}(P).
\]
Hence
\[
    \hat{\mathcal T}_{N,M}(P)
    =
    \mathcal R^M(P)+C_{N,M}(P),
\]
and adding and subtracting $\mathcal R^N(P)$ proves
\eqref{eq:finite_pi_perturbation}. Since $\Pi_i(P)\succeq0$, also
$C_{N,M}(P)\succeq0$.

Consequently,
\[
    \hat{\mathcal T}_{N,M}(P)
    \succeq
    \mathcal R^M(P)
    \succeq
    S_\nu,
\]
so the lower boundary is automatically preserved. By
\eqref{eq:strict_upper_slack},
\[
    D_M(P)
    =
    P_{\max}-\mathcal R^M(P)
    \succeq\eta I.
\]
Preservation of the upper boundary is therefore equivalent to
\(
    C_{N,M}(P)\preceq D_M(P),
\)
or, equivalently,
\[
    \lambda_{\max}\!\left(
        D_M(P)^{-1/2}
        C_{N,M}(P)
        D_M(P)^{-1/2}
    \right)
    \le1.
\]
This holds for every $P\in\mathcal P_\nu$ if and only if
$\rho_{N,M}\le1$.

Finally, the objective in \eqref{eq:rho_NM} is continuous in $P$:
all finite Riccati iterates, greedy gains, and $C_{N,M}(P)$ are
continuous, while $D_M(P)\succeq\eta I$. Since $\mathcal P_\nu$ is
compact, the maximum is finite and attained.
\end{proof}

The perturbation consists of two terms: the difference
between the evaluation and planning horizons $M$ and $N$,
and the additional cost of evaluating the fixed RHC gain
$K_N(P)$ instead of reoptimizing at each Bellman step.
The latter is captured by $C_{N,M}(P)\succeq0$.

The quantities in Lemma~\ref{lem:robust_pi_invariance} are
finite-horizon objects. For every fixed $P\in\mathcal P_\nu$,
evaluation of the objective in \eqref{eq:rho_NM} requires only the
finite Riccati iterates
\[
    \mathcal R^i(P),
    \qquad
    i=0,\ldots,\max\{N-1,M\},
\]
together with the corresponding greedy gains and finite-dimensional
matrix operations. The outer maximization in \eqref{eq:rho_NM} is
therefore a finite-dimensional optimization over the compact set
$\mathcal P_\nu$; no DARE solution is required.

Similarly, define the exact worst-case perturbation from the $N$-step
planner by
\begin{equation}
\label{eq:e_NM_exact}
    e_{N,M}^{\max}
    :=
    \max_{P\in\mathcal P_\nu}
    d_T\!\left(
        \hat{\mathcal T}_{N,M}(P),
        \mathcal R^N(P)
    \right).
\end{equation}
By continuity and compactness, $e_{N,M}^{\max}$ is finite and attained.
Thus $\rho_{N,M}\le1$ is an exact finite-dimensional invariance
certificate, while $e_{N,M}^{\max}$ is the corresponding exact
worst-case perturbation from the nominal planner.

A conservative bound on $e_{N,M}^{\max}$ follows directly from the
margin certificate. Let
\[
    P_M^-:=\mathcal R^M(S_\nu),
    \quad
    \kappa_M
    :=
    \lambda_{\max}\!\left(
        (P_M^-)^{-1/2}
        P_{\max}
        (P_M^-)^{-1/2}
    \right),
\]
and define the horizon-mismatch term
\[
    \Delta_{N,M}
    :=
    \max_{P\in\mathcal P_\nu}
    d_T\!\left(
        \mathcal R^M(P),
        \mathcal R^N(P)
    \right).
\]
Since
\[
    P_{\max}
    \preceq
    \kappa_M P_M^-
    \preceq
    \kappa_M\mathcal R^M(P),
\]
we have
\[
    D_M(P)
    \preceq
    (\kappa_M-1)\mathcal R^M(P).
\]
Hence, whenever $\rho_{N,M}\le1$,
\[
    C_{N,M}(P)
    \preceq
    \rho_{N,M}D_M(P)
    \preceq
    \rho_{N,M}(\kappa_M-1)\mathcal R^M(P),
\]
and therefore
\[
    \mathcal R^M(P)
    \preceq
    \hat{\mathcal T}_{N,M}(P)
    \preceq
    \bigl(1+\rho_{N,M}(\kappa_M-1)\bigr)
    \mathcal R^M(P).
\]
It follows that
\[
    d_T\!\left(
        \hat{\mathcal T}_{N,M}(P),
        \mathcal R^M(P)
    \right)
    \le
    \log\!\left(
        1+\rho_{N,M}(\kappa_M-1)
    \right).
\]
Thus
\begin{equation}
\label{eq:e_NM_bound}
\boxed{
    e_{N,M}^{\max}
    \le
    \bar e_{N,M},
}
\end{equation}
where \(    \bar e_{N,M}
    :=
    \Delta_{N,M}
    +
    \log\!\left(
        1+\rho_{N,M}(\kappa_M-1)
    \right)
    \).

From the $\nu$-step Riccati contraction established above,
\[
    d_T(\mathcal R^\nu(P_1),\mathcal R^\nu(P_2))
    \le
    \bar\gamma_\nu d_T(P_1,P_2),
    \qquad
    \bar\gamma_\nu<1.
\]
Since $\mathcal R$ is order preserving and subhomogeneous,
$\mathcal R(\zeta P)\preceq\zeta\mathcal R(P)$ for $\zeta\ge1$,
it is Thompson-nonexpansive. Indeed, if
$\delta=d_T(P,P')$, then
\(
    P\preceq e^\delta P',
\) and
\(
    P'\preceq e^\delta P,
\)
and therefore
\[
    \mathcal R(P)
    \preceq
    \mathcal R(e^\delta P')
    \preceq
    e^\delta\mathcal R(P'),
\]
with the reverse inequality obtained symmetrically. Hence
\[
    d_T(\mathcal R(P),\mathcal R(P'))
    \le
    d_T(P,P'),
\]
and the same holds for any finite composition $\mathcal R^r$.
Hence, for the $N$-step planner with $N\ge\nu$ and
$r=N-\lfloor N/\nu\rfloor\nu$,
\[
    d_T(\mathcal R^N(P_1),\mathcal R^N(P_2))
    \le
    \gamma_N d_T(P_1,P_2),
    \quad
    \gamma_N
    :=
    \bar\gamma_\nu^{\lfloor N/\nu\rfloor}
    <1.
\]

We can now quantify the effect of truncated policy evaluation on the
contractive $N$-step planner.

\begin{theorem}[ISS-type contraction of truncated PI]
\label{thm:iss_truncated_pi}
Let $N\ge\nu$, $M\ge1$ and
\(
    \hat{\mathcal T}_{N,M}(\mathcal P_\nu)
    \subseteq
    \mathcal P_\nu.
\)

For any $P_0\in\mathcal P_\nu$, consider
\(
    P_{k+1}
    =
    \hat{\mathcal T}_{N,M}(P_k),
\)
and define
\[
    d_k:=d_T(P_k,P^\star),
    \qquad
    e_k:=
    d_T(P_{k+1},\mathcal R^N(P_k)).
\]
Then $P_k\in\mathcal P_\nu$ for all $k\ge0$, and
\begin{equation}
\label{eq:iss_pi}
\boxed{
    d_{k+1}
    \le
    \gamma_N d_k+e_k,
    \qquad
    e_k
    \le
    e_{N,M}^{\max}
    \le
    \bar e_{N,M}.
}
\end{equation}
Consequently,
\begin{equation}
\label{eq:iss_pi_iterated}
    d_k
    \le
    \gamma_N^k d_0
    +
    \sum_{\ell=0}^{k-1}
    \gamma_N^{k-1-\ell}e_\ell,
\end{equation}
and, in particular,
\begin{equation}
\label{eq:iss_ultimate_bound}
\boxed{
    \limsup_{k\to\infty}d_k
    \le
    \frac{e_{N,M}^{\max}}{1-\gamma_N}
    \le
    \frac{\bar e_{N,M}}{1-\gamma_N}.
} \vspace{1em}
\end{equation}
\end{theorem}

\begin{proof}
By Lemma~\ref{lem:robust_pi_invariance}, $P_k\in\mathcal P_\nu$ for all $k\ge0$ whenever
$P_0\in\mathcal P_\nu$.
Moreover, $\mathcal R^N$ is a strict Thompson-metric contraction on
$\mathcal P_\nu$ with factor $\gamma_N<1$, and $P^\star$ is its unique
fixed point. Thus
\[
\begin{aligned}
d_{k+1}
&=
d_T(P_{k+1},P^\star) \\
&\le
d_T(P_{k+1},\mathcal R^N(P_k))
+
d_T(\mathcal R^N(P_k),\mathcal R^N(P^\star)) \\
&\le
e_k+\gamma_N d_k.
\end{aligned}
\]
By \eqref{eq:e_NM_exact} and \eqref{eq:e_NM_bound},
\(
    e_k
    \le
    e_{N,M}^{\max}
    \le
    \bar e_{N,M},
\)
which proves \eqref{eq:iss_pi}.

Repeated substitution gives \eqref{eq:iss_pi_iterated}. Using
$e_\ell\le e_{N,M}^{\max}$ yields
\[
    d_k
    \le
    \gamma_N^k d_0
    +
    \frac{1-\gamma_N^k}{1-\gamma_N}
    e_{N,M}^{\max}.
\]
Taking the limit superior gives
\[
    \limsup_{k\to\infty}d_k
    \le
    \frac{e_{N,M}^{\max}}{1-\gamma_N},
\]
and \eqref{eq:e_NM_bound} gives the second inequality in
\eqref{eq:iss_ultimate_bound}.
\end{proof}

Neither the truncated PI update nor the robust-invariance and
ultimate-bound certificates require solving the DARE. The interval
\(
    \mathcal P_\nu=[S_\nu,P_{\max}]
\)
is constructed using finite Riccati iterations and Lyapunov equations.
Moreover,
\[
    d_T(P,P^\star)
    \le
    \frac{
        d_T(P,\mathcal R^N(P))
    }{
        1-\gamma_N
    },
    \qquad
    P\in\mathcal P_\nu,
\]
provides a computable a-posteriori fixed-point error estimate. Indeed,
since $P^\star=\mathcal R^N(P^\star)$,
\[
\begin{aligned}
d_T(P,P^\star)
&\le
d_T(P,\mathcal R^N(P))
+
d_T(\mathcal R^N(P),\mathcal R^N(P^\star)) \\
&\le
d_T(P,\mathcal R^N(P))
+
\gamma_N d_T(P,P^\star),
\end{aligned}
\]
and rearranging gives the stated bound. Thus the computable
fixed-point residual provides an a-posteriori certificate of the
distance to $P^\star$, which remains an analysis object rather than an
algorithmic input.

The ISS estimate provides a quantitative bound on the asymptotic
distance to $P^\star$, but does not by itself determine whether the
truncated-PI iteration converges exactly to $P^\star$. We therefore
examine the local structure of the truncated-PI map and subsequently
combine it with the Thompson-metric contraction estimate.

\begin{theorem}[Local structure of truncated PI]
\label{thm:local_asymptotic_pi}
For any finite horizons $N\ge1$ and $M\ge1$,
\[
    \hat{\mathcal T}_{N,M}(P^\star)=P^\star.
\]
Moreover, the Fr\'echet derivative of the truncated-PI map at
$P^\star$ is
\begin{equation}
\label{eq:local_pi_derivative}
    D\hat{\mathcal T}_{N,M}(P^\star)[\Xi]
    =
    ((A^\star)^M)^\top \Xi(A^\star)^M,
\end{equation}
where
\(
    K^\star:=G(P^\star),
\)
\(
    A^\star:=A+BK^\star.
\)
Consequently,
\[
    \rho\!\left(
        D\hat{\mathcal T}_{N,M}(P^\star)
    \right)
    =
    \rho(A^\star)^{2M}
    <1.
\]
In addition, the perturbation term in
Lemma~\ref{lem:robust_pi_invariance} satisfies
\[
    C_{N,M}(P)
    =
    O\!\left(\|P-P^\star\|_F^2\right) \qquad \text{as } P\to P^\star.
\]

\end{theorem}

\begin{proof}
Since $\mathcal R(P^\star)=P^\star$,
\[
    K_N(P^\star)
    =
    G(\mathcal R^{N-1}(P^\star))
    =
    K^\star.
\]
Hence truncated evaluation of $K^\star$ initialized with $P^\star$
returns $P^\star$, so
\(
    \hat{\mathcal T}_{N,M}(P^\star)=P^\star.
\)

By Lemma~\ref{lem:robust_pi_invariance},
\(
    \hat{\mathcal T}_{N,M}(P)
    =
    \mathcal R^M(P)+C_{N,M}(P).
\)
At $P=P^\star$, every greedy gain appearing in $\Pi_i(P)$ equals
$K^\star$. Moreover, the corresponding gain mismatch is
$O(\|P-P^\star\|_F)$, while each $\Pi_i(P)$ is quadratic in this
mismatch. Hence
\(
    C_{N,M}(P)
    =
    O\!\left(\|P-P^\star\|^2\right),
\)
and therefore
\(
    DC_{N,M}(P^\star)[\Xi]=0.
\)
It follows that
\(
    D\hat{\mathcal T}_{N,M}(P^\star)[\Xi]
    =
    D\mathcal R^M(P^\star)[\Xi].
\)

The derivative of the Riccati map at its fixed point is
\[
    D\mathcal R(P^\star)[\Xi]
    =
    (A^\star)^\top \Xi A^\star.
\]
Applying the chain rule $M$ times gives
\[
    D\mathcal R^M(P^\star)[\Xi]
    =
    ((A^\star)^M)^\top \Xi(A^\star)^M,
\]
which proves \eqref{eq:local_pi_derivative}. Since $A^\star$ is Schur,
\[
    \rho\!\left(
        D\hat{\mathcal T}_{N,M}(P^\star)
    \right)
    =
    \rho(A^\star)^{2M}
    <1.
\]
\end{proof}

The first-order behavior in \eqref{eq:local_pi_derivative} depends on
the evaluation horizon $M$, but not on the prediction horizon $N$.
In particular, the spectral radius of the linearized map is
$\rho(A^\star)^{2M}$, which decreases as $M$ increases. The prediction
horizon $N$ therefore affects the truncated-PI map only through
higher-order terms around $P^\star$.
We next combine this local second-order structure with the
Thompson-metric contraction estimate. This yields a computable regional
certificate for exact convergence to $P^\star$, rather than merely an
ISS neighborhood.

\begin{corollary}[Computable regional asymptotic convergence]
\label{cor:computable_regional_convergence}
Let $N\ge\nu$ and $M\ge\nu$, and suppose
\(
    \rho_{N,M}\le1.
\)
Define
\[
    r_M(P)
    :=
    d_T(P,\mathcal R^M(P)),
    \;
    e_M(P)
    :=
    d_T(\hat{\mathcal T}_{N,M}(P),\mathcal R^M(P)).
\]
For $\delta>0$, define
\[
    \Omega_\delta
    :=
    \{P\in\mathcal P_\nu:r_M(P)\le\delta\}.
\]
For the nontrivial case $P\neq P^\star$, for which $r_M(P)>0$, define
\[
    \omega_M(\delta)
    :=
    \sup_{P\in\Omega_\delta\setminus\{P^\star\}}
    \frac{e_M(P)}{r_M(P)}.
\]
Then
\(
    \omega_M(\delta)\to0
\)
as
\(
    \delta\downarrow0.
\)
If
\begin{equation}
\label{eq:local_capture_test}
    \omega_M(\delta)
    <
    \frac{1-\gamma_M}{1+\gamma_M},
\end{equation}
then $\Omega_\delta$ is forward invariant and
\[
    r_M(P_{k+1})
    \le
    q_\delta r_M(P_k),
    \qquad
    q_\delta
    :=
    \gamma_M+(1+\gamma_M)\omega_M(\delta)<1.
\]
Hence every trajectory entering $\Omega_\delta$ converges
to $P^\star$.

Moreover, if
\begin{equation}
\label{eq:iss_capture_condition}
    \frac{e_{N,M}^{\max}}{1-\gamma_N}
    <
    \frac{\delta}{1+\gamma_M},
\end{equation}
then the approximate PI initialized in
$\mathcal P_\nu$ converges to $P^\star$.
\end{corollary}

\begin{proof}
Since $\mathcal R^M$ is a $\gamma_M$-contraction,
\[
    (1-\gamma_M)d_T(P,P^\star)
    \le
    r_M(P)
    \le
    (1+\gamma_M)d_T(P,P^\star).
\]
By Theorem~\ref{thm:local_asymptotic_pi},
\[
    C_{N,M}(P)
    =
    O\!\left(\|P-P^\star\|_F^2\right).
\]
By Proposition~1.1 of \cite{nussbaum1994finsler}, the Thompson metric
is locally equivalent to the matrix norm on the interior of the cone.
Since $\mathbb S^n$ is finite-dimensional, it is therefore locally
equivalent to the Frobenius norm around $P^\star$. Using
\[
    \hat{\mathcal T}_{N,M}(P)-\mathcal R^M(P)
    =
    C_{N,M}(P),
\]
and noting that both $\mathcal R^M(P)$ and
$\hat{\mathcal T}_{N,M}(P)$ converge to $P^\star$ as
$P\to P^\star$, we obtain
\[
    e_M(P)
    =
    O\!\left(d_T(P,P^\star)^2\right).
\]

Moreover, the lower residual bound gives
\[
    d_T(P,P^\star)
    \le
    \frac{r_M(P)}{1-\gamma_M}.
\]
Hence, for $P\in\Omega_\delta$,
\[
    d_T(P,P^\star)
    \le
    \frac{\delta}{1-\gamma_M},
\]
so $\Omega_\delta$ shrinks to $P^\star$ as $\delta\downarrow0$.
Combining the preceding estimates yields
\[
    e_M(P)=O\!\left(r_M(P)^2\right),
\]
and therefore
\[
    \omega_M(\delta)\to0
    \qquad
    \text{as }\delta\downarrow0.
\]

Using the triangle inequality and contraction of $\mathcal R^M$,
\[
\begin{aligned}
    r_M(P_{k+1})
    &=
    d_T(P_{k+1},\mathcal R^M(P_{k+1}))\\
    &\le
    e_M(P_k)
    +
    \gamma_M d_T(P_k,P_{k+1})\\
    &\le
    \gamma_M r_M(P_k)
    +(1+\gamma_M)e_M(P_k).
\end{aligned}
\]
Therefore, for $P_k\in\Omega_\delta$,
\[
    r_M(P_{k+1})
    \le
    q_\delta r_M(P_k).
\]
Condition \eqref{eq:local_capture_test} gives $q_\delta<1$, so
$\Omega_\delta$ is forward invariant and $r_M(P_k)\to0$
geometrically. The lower bound above then implies
\(
    d_T(P_k,P^\star)\to0,
\)
and hence $P_k\to P^\star$.

Finally, Theorem~\ref{thm:iss_truncated_pi} and
\eqref{eq:iss_capture_condition} give
\[
    \limsup_{k\to\infty}d_T(P_k,P^\star)
    <
    \frac{\delta}{1+\gamma_M}.
\]
Thus, for all sufficiently large $k$,
\[
    r_M(P_k)
    \le
    (1+\gamma_M)d_T(P_k,P^\star)
    <
    \delta.
\]
The trajectory therefore enters $\Omega_\delta$ in finite time and
subsequently converges geometrically to $P^\star$.
\end{proof}

Increasing $N$ strengthens the contraction of the nominal Riccati
planner through $\gamma_N$, at the expense of a longer online
look-ahead optimization. The evaluation horizon $M$ affects the
finite-rollout perturbation, the invariance condition
$\rho_{N,M}\le1$, and the local convergence mechanism without
increasing the online planning horizon. For a fixed $N$, regional
asymptotic convergence can therefore be certified by finding a finite
$M$ and $\delta>0$ such that
\[
    \rho_{N,M}\le1,
    \qquad
    \omega_M(\delta)
    <
    \frac{1-\gamma_M}{1+\gamma_M},
    \qquad
    \frac{e_{N,M}^{\max}}{1-\gamma_N}
    <
    \frac{\delta}{1+\gamma_M}.
\]
If no such certificate is obtained for the chosen $N$, the prediction
horizon can be increased and the test repeated.

\begin{remark}[Online diagnostics]
\label{rem:online_diagnostics}
The quantities $\rho_{N,M}$, $e_{N,M}^{\max}$, and $\gamma_N$ are
uniform worst-case quantities over $\mathcal P_\nu$. Along a realized
trajectory,
\[
    \rho_k
    :=
    \lambda_{\max}\!\left(
        D_M(P_k)^{-1/2}
        C_{N,M}(P_k)
        D_M(P_k)^{-1/2}
    \right)
    \le
    \rho_{N,M},
\]
and therefore $\rho_k$ quantifies how much of the available upper order
margin is used by the realized truncated-PI correction.
More importantly for convergence,
\[
    e_k
    =
    d_T(P_{k+1},\mathcal R^N(P_k))
\]
is the realized perturbation from the nominal $N$-step Riccati planner.
Whenever $P_k\neq P^\star$, define
\[
    \gamma_k
    :=
    \frac{
        d_T(\mathcal R^N(P_k),P^\star)
    }{
        d_T(P_k,P^\star)
    }.
\]
Since $\mathcal R^N$ is a $\gamma_N$-contraction,
\(
    0\le\gamma_k\le\gamma_N<1,
\)
and the triangle inequality gives the trajectory-wise estimate
\[
    d_{k+1}
    \le
    \gamma_k d_k+e_k.
\]
Hence the gaps between the uniform and realized quantities quantify the
conservatism of the regional certificates. Only $\gamma_k$ requires
knowledge of $P^\star$.
\end{remark}

\subsection{Closed-Loop Stability under Policy Updates}
The convergence analysis above concerns the terminal cost iteration. For online deployment, we additionally need to ensure that stability is preserved when the induced RHC policy is updated. First, we establish sufficient conditions for closed-loop stability under policy updates. We then show how reduction of the step size and dwell-time safeguards can enforce these conditions throughout closed-loop updates.

\begin{theorem}[Stability under policy updates]
\label{thm:truncated_pi_stability}
Let
\(
    K_k:=K_N(P_k),
\) and
\(
    A_k:=A+BK_k.
\)
Assume that $A_0$ is Schur. For each $k$, let $X_k\succ0$ solve
\[
    X_k-A_k^\top X_kA_k=I.
\]
If
\begin{equation}
\label{eq:lyap_x}
    X_k-A_{k+1}^\top X_kA_{k+1}\succ0,
    \qquad k\ge0,
\end{equation}
then every updated policy $K_k$ is stabilizing.

Moreover, let $\bar P_k\succ0$ denote the infinite-horizon value matrix
associated with $K_k$,
\[
    \bar P_k
    =
    Q+K_k^\top RK_k+A_k^\top\bar P_kA_k,
\]
and let $K_k$ be applied for $\tau_k$ time steps before activating
$K_{k+1}$. If there exists $\bar\vartheta\in(0,1)$ such that
\begin{equation}
\label{eq:lyap_w}
    (A_k^{\tau_k})^\top
    \bar P_{k+1}
    A_k^{\tau_k}
    \preceq
    \bar\vartheta\bar P_k,
    \qquad k\ge0,
\end{equation}
then the Lyapunov value decreases geometrically at the policy-update
instants. If, in addition, $\sup_k\tau_k<\infty$, the resulting
time-varying closed loop is exponentially stable in time $t$.
\end{theorem}

\begin{proof}
Suppose $A_k$ is Schur. Then $X_k\succ0$ is well defined, and
\eqref{eq:lyap_x} makes $X_k$ a strict Lyapunov matrix for $A_{k+1}$.
Hence $A_{k+1}$ is Schur, and the first claim follows by induction.

Let $t_{k+1}=t_k+\tau_k$ denote the policy-update times and define
$V_k(x):=x^\top\bar P_kx$. Since
$x_{t_{k+1}}=A_k^{\tau_k}x_{t_k}$, \eqref{eq:lyap_w} gives
\[
    V_{k+1}(x_{t_{k+1}})
    \le
    \bar\vartheta V_k(x_{t_k}),
\]
and hence
\[
    V_k(x_{t_k})
    \le
    \bar\vartheta^kV_0(x_0).
\]
Since $\bar P_k\succeq P^\star\succeq S_\nu\succ0$,
\[
    \|x_{t_k}\|^2
    \le
    \frac{\lambda_{\max}(\bar P_0)}
         {\lambda_{\min}(S_\nu)}
    \bar\vartheta^k\|x_0\|^2.
\]
Between updates, $V_k$ is nonincreasing. Hence, if
$\sup_k\tau_k<\infty$, geometric decay in the update index $k$
implies exponential decay in plant time $t$.
\end{proof}

Condition \eqref{eq:lyap_x} is sufficient for stability of the updated policy and can always be enforced by sufficiently damping the terminal-cost update.

\begin{remark}[Safeguarding RHC policy updates]
\label{rem:pi_safeguard}
For the nominal truncated-PI candidate
\[
    \widetilde P_{k+1}
    =
    \hat{\mathcal T}_{N,M}(P_k),
\]
condition~\eqref{eq:lyap_x} can be checked before the corresponding
policy is activated. If it fails, consider the damped update
\[
    P_{k+1}(\alpha)
    =
    (1-\alpha)P_k+\alpha\widetilde P_{k+1},
    \qquad
    \alpha\in(0,1].
\]
Since
\(
    K_N(P_{k+1}(\alpha))
    \to K_k
    \)
    as 
    \(\alpha\downarrow0,
\)
and
\[
    X_k-A_k^\top X_kA_k=I\succ0,
\]
continuity guarantees that~\eqref{eq:lyap_x} holds for all sufficiently
small $\alpha>0$. Thus stability of the next policy can always be
recovered by a line search.

Likewise, if~\eqref{eq:lyap_w} fails at the intended policy-update
time, the current policy can be retained longer. Since $A_k$ is Schur,
\[
    (A_k^\tau)^\top\bar P_{k+1}A_k^\tau\to0
    \qquad(\tau\to\infty),
\]
so a finite dwell time satisfying~\eqref{eq:lyap_w} always exists.
Thus damping safeguards stability of the updated RHC policy,
while a sufficient dwell time safeguards switching between policies.
\end{remark}

\begin{remark}[Effect of damping on the contraction estimate]
\label{rem:damping_contraction}
When damping is used to safeguard online RHC policy updates,
it modifies the nominal truncated-PI iteration, but its effect enters
the preceding contraction analysis as an additional perturbation. With
\[
    P_{k+1}
    =
    (1-\alpha_k)P_k+\alpha_k\widetilde P_{k+1},
\]
the triangle inequality gives
\[
\begin{aligned}
d_T(P_{k+1},\mathcal R^N(P_k))
&\le
d_T(P_{k+1},\widetilde P_{k+1})
+
d_T(\widetilde P_{k+1},\mathcal R^N(P_k))\\
&=:e_k^{\rm safe}+\widetilde e_k.
\end{aligned}
\]
As $\mathcal P_\nu$ is convex, damping preserves
$\mathcal P_\nu$ whenever
$P_k,\widetilde P_{k+1}\in\mathcal P_\nu$.
Hence the ISS-type estimate applies with the nominal
perturbation $\widetilde e_k$ replaced by
$\widetilde e_k+e_k^{\rm safe}$.
When damping is inactive, $e_k^{\rm safe}=0$ and the nominal contraction
estimate is recovered exactly.

Since \eqref{eq:lyap_x} holds with strict margin at $P^\star$,
continuity ensures that the full step $\alpha_k=1$ is admissible
in a neighborhood of $P^\star$. Since the nominal contraction and full-step admissibility hold simultaneously in a sufficiently small neighborhood of $P^\star$, the safeguarded RHC iteration is locally convergent. More generally, local convergence under damping is retained if $\sum_k\alpha_k=\infty$. Consequently, any safeguarded sequence converging to $P^\star$ eventually accepts full steps. The numerical examples in Section~\ref{sec:numerical_results} illustrate this behavior.
\end{remark}

\section{Model Predictive Control}
\label{sec:mpc}
Here, we extend the PI framework to the MPC setting. In the
constrained case, however, the value function cannot generally
be represented by a single matrix, e.g., under affine constraints
it is piecewise quadratic. The RHC analysis provides candidate
terminal cost updates, while terminal pair compatibility and
safe switching conditions are used to preserve recursive
feasibility and the MPC Lyapunov decrease.

\subsection{Problem Formulation}
The MPC policy is (implicitly) obtained by solving a finite-horizon optimal control problem with a terminal cost and a terminal constraint set. At PI step \(k\) and time step $t$, we solve the following optimal control problem:
\begin{equation}
\label{eq:mpc}
    \begin{aligned}
        \min_{\mathbf{u}_{0 \mid t}, \ldots, \mathbf{u}_{N-1 \mid t}} & \sum_{i=0}^{N-1} \mathbf{x}_{i \mid t}^\top Q \mathbf{x}_{i \mid t} + \mathbf{u}_{i \mid t}^\top R \mathbf{u}_{i \mid t} + \mathbf{x}_{N \mid t}^\top P_k \mathbf{x}_{N \mid t} \\
        \text{s.t. } \mathbf{x}_{i+1 \mid t} & = A \mathbf{x}_{i \mid t} + B \mathbf{u}_{i \mid t}, \ i = 0, \ldots, N-1 \\
        & \mathbf{u}_{i \mid t} \in \mathbb{U}, \qquad i = 0, \ldots, N-1 \\
        & \mathbf{x}_{i \mid t} \in \mathbb{X}, \qquad i = 0, \ldots, N-1 \\
        & \mathbf{x}_{N \mid t} \in \mathbb{X}_{f,k},
    \end{aligned}
\end{equation}
where $\mathbf{x}_{i \mid t}$ is the prediction of $\mathbf{x}_{t+i}$ at time $t$. 
The optimal finite-horizon value of \eqref{eq:mpc} is denoted by $J_k^N(x)$, and depends on the terminal pair $(P_k,\mathbb X_{f,k})$.
The resulting MPC policy is denoted by $\pi_k^N(x)$.
The truncated RHC update provides a candidate terminal matrix
\(\widetilde P_{k+1}\). For constrained MPC, however, the terminal
matrix and terminal set must form a compatible terminal pair, as defined next.

For a terminal matrix \(P_k\), define
\[
K_{f,k}:=G(P_k),
\quad
A_{f,k}:=A+BK_{f,k},
\quad
V_{f,k}(x):=x^\top P_kx.
\]
We call \((P_k,\mathbb X_{f,k})\) a compatible terminal pair if
\begin{equation}
\label{eq:mpc_terminal_pair}
\begin{aligned}
\mathbb X_{f,k} &\subseteq \mathbb X, && \text{where } 0\in\operatorname{int}(\mathbb X_{f,k})\\
K_{f,k}x &\in\mathbb U,
&& x\in\mathbb X_{f,k},\\
A_{f,k}x &\in\mathbb X_{f,k},
&& x\in\mathbb X_{f,k},
\end{aligned}
\end{equation}
and
\begin{equation}
\label{eq:mpc_terminal_clf}
V_{f,k}(A_{f,k}x)-V_{f,k}(x)
\le
-\ell(x,K_{f,k}x),
\qquad
x\in\mathbb X_{f,k}.
\end{equation}
For \(K_{f,k}=G(P_k)\), the terminal decrease condition is guaranteed by
\begin{equation}
\label{eq:mpc_terminal_riccati}
\boxed{
\mathcal R(P_k)\preceq P_k,
}
\end{equation}
since
\(
\mathcal R(P_k)
=
Q+K_{f,k}^\top RK_{f,k}
+A_{f,k}^\top P_kA_{f,k}.
\)
Consequently, a strict Riccati subsolution cannot be used as a
quadratic MPC terminal cost: if $P\prec\mathcal R(P)$, then
\[
Q+K^\top RK+A_K^\top PA_K
\succeq
\mathcal R(P)
\succ P,
\qquad \forall K.
\]
We therefore initialize from a strict Riccati supersolution \eqref{eq:supersolution}.

\begin{remark}[Terminal set construction]
A standard choice is a quadratic sublevel set of the terminal cost:
\begin{equation}
\label{eq:mpc_quadratic_terminal_set}
\mathbb X_{f,k}
=
\{x:x^\top P_kx\le c_{f,k}\}.
\end{equation}
If \(\mathcal R(P_k)\preceq P_k\), such a set is automatically
positively invariant under \(K_{f,k}\); \(c_{f,k}\) can then be chosen
to satisfy the state and input constraints. Alternatively, more general terminal sets, including polytopic or maximal positively invariant sets, may be used to reduce conservatism.
\end{remark}

\subsection{Value Function Update from Constrained Data}
\label{sec:mpc_value_update}
The constrained MPC value function is generally non-quadratic, so the truncated-PI update cannot be applied directly. 
However, as the unconstrained RHC policy is locally admissible in a neighborhood of the origin, the unconstrained rollout can be used to construct a candidate update for the terminal cost.
In particular, consider an \(M\)-step rollout window for which the
unconstrained RHC optimizer is feasible for the constrained MPC problem
at every step. Over such a window, the two policies coincide and
\(
u_{t+\tau}=K_N(P_k)x_{t+\tau},
\)
for \(
\tau=0,1,\ldots,M-1,
\)
and the \(M\)-step value sample satisfies \eqref{eq:sample_value_target}.
Hence the same sample-based quadratic fit of
Section~\ref{sec:sample_based_value_fit} can be used to construct the
candidate update \eqref{eq:truncated_pi_update}.
Samples for which the constrained and unconstrained optimizers differ
generally bias a single quadratic fit; this effect is illustrated in
Section~\ref{sec:numerical_results}.

\subsection{Safe Update of Terminal Ingredients}
\label{sec:mpc_safe_update}
Although the RHC analysis characterizes the candidate terminal cost update, constrained MPC additionally requires compatibility of the terminal set and a feasible transition between consecutive terminal pairs.
To safely update the terminal ingredients, we require the terminal set update condition
\begin{equation}
\label{eq:mpc_set_switch}
\boxed{
A_{f,k}\mathbb X_{f,k}
\subseteq
\mathbb X_{f,k+1}.
}
\end{equation}
Moreover, to preserve the Lyapunov decrease across the update, we
require
\begin{equation}
\label{eq:mpc_cross_update_set}
V_{f,k+1}(A_{f,k}x)-V_{f,k}(x)
\le
-\ell(x,K_{f,k}x),
\quad
x\in\mathbb X_{f,k}.
\end{equation}
For the quadratic terminal costs considered here, a sufficient
condition for \eqref{eq:mpc_cross_update_set} is
\begin{equation}
\label{eq:mpc_cross_update}
\boxed{
A_{f,k}^\top P_{k+1}A_{f,k}
+
Q+K_{f,k}^\top RK_{f,k}
\preceq
P_k.
}
\end{equation}

\begin{theorem}[Safe MPC policy iteration]
\label{thm:safe_mpc_pi}
Let
\((P_k,\mathbb X_{f,k})\) and
\((P_{k+1},\mathbb X_{f,k+1})\)
be compatible terminal pairs satisfying
\eqref{eq:mpc_set_switch} and
\eqref{eq:mpc_cross_update_set}.
If \eqref{eq:mpc} is feasible at \(x_t\), then, after applying
\(
u_t=\pi_k^N(x_t),
\)
the MPC problem with terminal pair
\((P_{k+1},\mathbb X_{f,k+1})\) is feasible at \(x_{t+1}\), and
\begin{equation}
\label{eq:mpc_cross_decrease}
\boxed{
J_{k+1}^N(x_{t+1})
\le J_k^N(x_t)-\ell(x_t,u_t).
}
\end{equation}
Consequently, any sequence of terminal pair updates satisfying these
conditions preserves recursive feasibility and the MPC Lyapunov
decrease and, hence, closed-loop stability.
\end{theorem}

\begin{proof}
Let \(z=\mathbf x_{N|t}^\star\in\mathbb X_{f,k}\) denote the terminal
state of an optimal solution at time \(t\). At \(t+1\), shift the
optimal input sequence and append \(K_{f,k}z\). By
\eqref{eq:mpc_set_switch},
\(
A_{f,k}z\in\mathbb X_{f,k+1},
\)
so the shifted sequence is feasible for the updated MPC problem.
Since $z\in\mathbb X_{f,k}$, \eqref{eq:mpc_cross_update_set} ensures
that the appended stage cost and new terminal cost do not exceed the
previous terminal cost. Hence the shifted candidate yields
\eqref{eq:mpc_cross_decrease}.

For $x_0$ in a neighborhood of the origin contained in
$\mathbb X_{f,0}$, terminal pair compatibility and the
Lyapunov decrease give
\[
x_t^\top P^\star x_t
\le J^N_{k(t)}(x_t)
\le J^N_0(x_0)
\le x_0^\top P_0x_0,
\]
which establishes Lyapunov stability. Moreover, summing the
decrease yields $\sum_{t=0}^\infty\ell(x_t,u_t)<\infty$,
and hence $x_t\to0$ by observability of $(A,Q^{1/2})$ and $R\succ0$.
\end{proof}

Subsequently, we propose two strategies to ensure that the established conditions for safe switching are satisfied. The first strategy is to reduce step size of the update until the conditions are satisfied for all states, while the second strategy is to wait until a state is reached that satisfies a state-dependent switching condition.

\begin{remark}[Damping]
\label{rem:mpc_damping}
If the candidate \(\widetilde P_{k+1}\) does not satisfy the terminal
or switching conditions, consider
\begin{equation}
\label{eq:mpc_damped_update}
P_{k+1}(\alpha)
=
(1-\alpha)P_k+\alpha\widetilde P_{k+1},
\qquad
\alpha\in(0,1].
\end{equation}
If \(P_k\) is a strict Riccati supersolution, then, by continuity,
sufficiently small positive \(\alpha\) preserves both
\(
\mathcal R(P_{k+1}(\alpha))
\preceq
P_{k+1}(\alpha)
\)
and \eqref{eq:mpc_cross_update}. Compatibility of the terminal set and the switching condition
\eqref{eq:mpc_set_switch} can be checked explicitly. If no acceptable
positive step is found, the current terminal pair can be retained until the next candidate update is feasible.
\end{remark}

\begin{remark}[Waiting]
\label{rem:mpc_waiting}
Instead of damping the candidate until an immediate switch is
certified, one may damp only as needed to obtain a compatible terminal
pair $(P_{k+1},\mathbb X_{f,k+1})$ and retain the current MPC pair
until a safe switching state is reached. The candidate pair may be activated at the first time $t$ after its generation satisfying
\begin{equation}
\label{eq:mpc_wait_switch}
x_t\in\mathbb X_{f,k+1},
\quad
V_{f,k+1}(x_t)\le\vartheta J_k^N(x_{t_k}),
\; \vartheta\in(0,1).
\end{equation}
The first condition ensures feasibility of the new MPC problem, while
compatibility of the new terminal pair gives
\[
    J_{k+1}^N(x_t)
    \le
    V_{f,k+1}(x_t)
    \le
    \vartheta J_k^N(x_{t_k}).
\]
Thus waiting guarantees a decrease between activation values, rather than the one-step cross-update decrease of Theorem~\ref{thm:safe_mpc_pi}.
Since the current MPC closed loop converges to the origin and $\mathbb X_{f,k+1}$ contains a neighborhood of the origin, such a switching time exists. Thus damping modifies the candidate to permit an immediate switch, whereas waiting retains the candidate and delays its activation.
\end{remark}

\begin{remark}[Ellipsoidal terminal sets]
Let the next terminal set be a quadratic sublevel set following \eqref{eq:mpc_quadratic_terminal_set}.
Then \eqref{eq:mpc_set_switch} is equivalent to
\[
A_{f,k}^\top P_{k+1}A_{f,k}
\preceq
\frac{c_{f,k+1}}{c_{f,k}}P_k.
\]
Hence, for ellipsoidal terminal sets, both the switching conditions reduce to finite-dimensional matrix tests.
\end{remark}

Damping permits an immediate safe switch, whereas waiting
delays activation until a suitable state is reached. The two
strategies are summarized in Algorithm~\ref{alg:constrained_policy_iteration}.
Here, $t_k$ denotes the activation time of policy $k$;
during waiting periods, the current policy remains active
and plant time continues to advance.

For constrained MPC, active constraints and switching
safeguards may prevent the local convergence argument
of Remark~\ref{rem:damping_contraction} from applying,
even when the corresponding unconstrained RHC iteration
converges.

\begin{algorithm}[t]
\caption{Online terminal cost policy iteration for RHC and MPC}
\label{alg:constrained_policy_iteration}
\begin{algorithmic}[1]
    \STATE Given $\mathtt{mode}\in
    \{\mathrm{RHC},\mathrm{MPC\text{-}damp},
    \mathrm{MPC\text{-}wait}\}$, $P_0$, $N$ and $M$.
    For MPC, initialize a compatible pair
    $(P_0,\mathbb X_{f,0})$.
    Set $k=t=t_k=0$ and $\mathcal D_k=\varnothing$.

    \WHILE{the closed loop is running}
        \STATE Compute and apply the current RHC policy or, for MPC,
        solve \eqref{eq:mpc} and apply $u_t=\pi_k^N(x_t)$.

        \STATE Compute the $M$-step value sample according to
        \eqref{eq:sample_value_target}.
        For MPC, retain the sample only if the policy coincides with RHC over
        the complete $M$-step window.

        \IF{$\mathcal D_k$ contains sufficiently informative samples}
            \STATE Compute the candidate $\widetilde P_{k+1}$ from
            \eqref{eq:sample_based_lsq}.

            \IF{$\mathtt{mode}=\mathrm{RHC}$}
                \STATE Damp $\widetilde P_{k+1}$ following \eqref{eq:mpc_damped_update} until
                \eqref{eq:lyap_x} holds.
                \STATE Activate the new policy when
                $t-t_k\ge\tau_k$, with $\tau_k$
                satisfying \eqref{eq:lyap_w}.

            \ELSIF{$\mathtt{mode}=\mathrm{MPC\text{-}damp}$}
                \STATE Damp $\widetilde P_{k+1}$ following \eqref{eq:mpc_damped_update} and construct the
                corresponding terminal set such that \eqref{eq:mpc_terminal_riccati}, \eqref{eq:mpc_set_switch}, and
                \eqref{eq:mpc_cross_update} are satisfied.
                \STATE Update the terminal pair immediately.

            \ELSIF{$\mathtt{mode}=\mathrm{MPC\text{-}wait}$}
                \STATE Damp $\widetilde P_{k+1}$ following \eqref{eq:mpc_damped_update} until
                \eqref{eq:mpc_terminal_riccati} holds,
                and construct a compatible terminal set.
                \STATE Retain the current terminal pair until
                \eqref{eq:mpc_wait_switch} holds.
                \STATE Update the terminal pair.
            \ENDIF

            \STATE Upon activation, set $t_k\leftarrow t$, $k\leftarrow k+1$, and $\mathcal D_k\leftarrow\varnothing$.
        \ENDIF

        \STATE Set $t\leftarrow t+1$.
    \ENDWHILE
\end{algorithmic}
\end{algorithm}

\section{Numerical Example}
\label{sec:numerical_results}
We demonstrate the proposed RHC and MPC policy iteration schemes, which rely on the underlying contraction of the Riccati map, in a numerical example. 
We consider a sampled double-integrator system with
sampling time \(0.1\),
\[
A=
\begin{bmatrix}
1 & 0.1\\
0 & 1
\end{bmatrix},
\qquad
B=
\begin{bmatrix}
0.005\\
0.1
\end{bmatrix},
\]
and cost matrices
\(
Q=\operatorname{diag}(10,1) \) and \(R=1\).
For the RHC policy, we consider no input or state constraints, while for the MPC policy, we consider input and state constraints of the form $\mathbb{U} = \{ u : |u| \leq 5 \}$ and $\mathbb{X} = \{ x : |x_1| \leq 2, |x_2| \leq 1 \}$.

The corresponding stabilizing Riccati solution is
\[
P^\star=
\begin{bmatrix}
90.7756 & 31.6623\\
31.6623 & 27.6585
\end{bmatrix}.
\]

\begin{figure}[t]
    \centering
    \includegraphics[width=\linewidth]{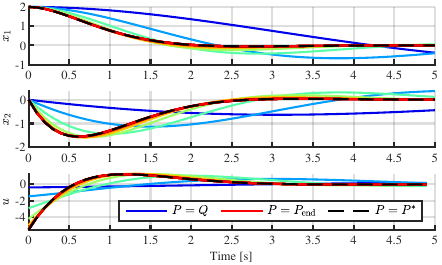}
    \caption{Receding-horizon control with episodic truncated policy
    iteration for \(N=M=2\). The terminal matrix is updated after
    each episode.}
    \label{fig:VFL}
\end{figure}

\subsection{Receding-Horizon Control}
For the RHC experiment, we choose \(N=M=2\) and the interval
\(
\mathcal P_1=[Q,5P^\star].
\)
Since \(Q\succ0\), \(S_1=Q\).
Moreover,
\(
\kappa_2\approx40.811, \
\Delta_{2,2}=0,
\)
where the latter holds exactly because \(M=N\).
The uniform quantities $\rho_{N,M}$ and $e_{N,M}^{\max}$ are estimated
offline by numerical maximization over $\mathcal P_\nu$, while
$\rho_k$, $e_k$, and $\gamma_k$ are evaluated along the realized
closed-loop trajectory.
Quantities requiring maximization over \(\mathcal P_1\) are approximated
numerically using MATLAB's \texttt{fmincon} from multiple initial
conditions. We obtain the numerical estimates
\(
\rho_{2,2}\approx0.04788,
\quad
e_{2,2}^{\max}\approx0.13426.
\)
Substituting the estimated $\rho_{2,2}$ into
\eqref{eq:e_NM_bound} gives
\[
\bar e_{2,2}
\approx
\log\!\left(1+\rho_{2,2}(\kappa_2-1)\right)
\approx1.0668.
\]
These estimates are consistent with the robust-invariance
condition and illustrate the conservatism of the
constructive perturbation bound.
For comparison, numerical maximization of the contraction factor of \(\mathcal R^2\) on \(\mathcal P_1\) gives
\(
\gamma_2^{\rm approx}\approx0.9088,
\)
whereas the analytical contraction factor is
\(
\gamma_2\approx0.9876.
\)

\begin{figure}[t]
    \centering
    \includegraphics[width=\linewidth]{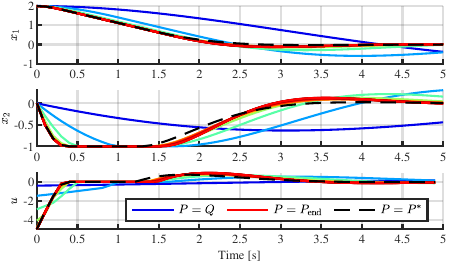}
    \caption{Biased MPC with constraints.}
    \label{fig:VFL_Biased}
\end{figure}

\begin{figure}[t]
    \centering
    \includegraphics[width=\linewidth]{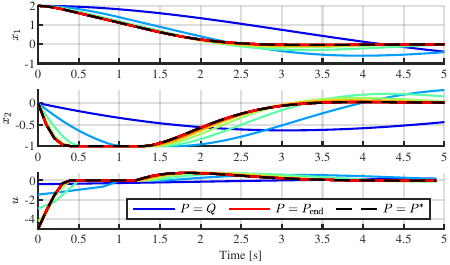}
    \caption{Unbiased MPC with constraints.}
    \label{fig:VFL_Constrained}
\end{figure}

We consider the two extremal initializations \(P_0=S_1\) and \(P_0=P_{\max}\).
Fig.~\ref{fig:VFL} shows the closed-loop trajectories generated by the
episodic truncated PI recursion for \(P_0 = S_1 = Q\).
We update the terminal cost after each episode for illustration purposes. Here, the closed-loop trajectories qualitatively converge from $P_0 = Q$ toward $P^\star$.
We compute the diagnostic statistics for the first fourteen updates, which are summarized in Table~\ref{tab:rhc_diagnostics}.
The realized contraction factors remain substantially below the
semi-global bound \(\gamma_2\), and the trajectory-wise invariance
factors satisfy \(\rho_k<1\) throughout both experiments. All iterates
remain in \(\mathcal P_1\), and all realized feedback matrices
\(A+BK_k\) are Schur.

\begin{figure}[t]
    \centering
    \includegraphics[width=\linewidth]{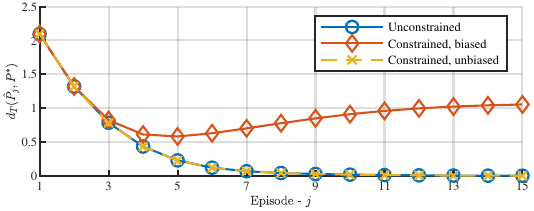}
    \caption{Thompson distance to $P^\star$ under terminal cost PI.}
    \label{fig:Raw_Thompson_Metric}
\end{figure}

The local exact-convergence mechanism of Corollary~\ref{cor:computable_regional_convergence} can also be observed numerically. In this example, $\frac{1-\gamma_2}{1+\gamma_2}\approx 6.26\times10^{-3}$. Numerical maximization over the residual sublevel set yields the estimate $\omega_2(0.0312) \approx 1.1024\times10^{-3} < (1-\gamma_2)/(1+\gamma_2)$, and hence the corresponding estimated contraction factor satisfies $q_\delta\approx0.9897<1$. Thus $\Omega_{0.0312}$ is a numerically identified candidate local attraction region for exact convergence. Both trajectories enter this region at iteration $k=8$ and exhibit numerical convergence to $P^\star$.

\begin{table}[h]
\centering
\caption{Diagnostics for the double-integrator RHC
example.}
\label{tab:rhc_diagnostics}
\begin{tabular}{lcc}
\toprule
 & \(P_0=S_1\) & \(P_0=P_{\max}\)\\
\midrule
\(d_T(P_0,P^\star)\)
    & \(3.4666\) & \(1.6094\)\\
\(d_T(P_{14},P^\star)\)
    & \(2.49\times10^{-3}\) & \(1.24\times10^{-3}\)\\
\(\max_k\gamma_k\)
    & \(0.6532\) & \(0.7827\)\\
\(\max_k e_k\)
    & \(6.72\times10^{-3}\) & \(6.63\times10^{-2}\)\\
\(\max_k\rho_k\)
    & \(6.45\times10^{-4}\) & \(4.79\times10^{-2}\)\\
\(\max_k\rho(A+BK_k)\)
    & \(0.9890\) & \(0.8735\)\\
\bottomrule
\end{tabular}
\vspace{-0.5em}
\end{table}

\subsection{Model Predictive Control}
\subsubsection{Biased Value Function Estimation}
We consider the same double-integrator system with state and input constraints. As discussed in Section~\ref{sec:mpc_value_update}, the constrained MPC value function cannot be represented as a single quadratic function, and a single quadratic fit based on constrained rollouts is generally biased. To recover the quadratic RHC update, we retain only rollout samples for which the unconstrained RHC optimizer is feasible for the constrained MPC problem.
Fig.~\ref{fig:VFL_Biased} and Fig.~\ref{fig:VFL_Constrained} compare the resulting closed-loop PI trajectories for biased and unbiased value function estimation. For the strict supersolution initialization \(P_0=P_{\max}\), the biased update can leave \(\mathbb S_{++}\), in which case the Thompson metric is not defined. Therefore, for this identification illustration, all recursions are initialized with \( P_0=Q \),
and we consider the closed-loop safety and stability guarantees in a subsequent numerical example.
Fig.~\ref{fig:Raw_Thompson_Metric} shows the Thompson metric between the estimated terminal matrix and $P^\star$.
The numerical results show that the unbiased estimator approaches the optimal LQR terminal matrix, whereas the biased estimator converges to a suboptimal solution.

\subsubsection{Closed-Loop Behavior with Safe Updates}
\begin{figure}[t]
    \centering
    \includegraphics[width=\linewidth]{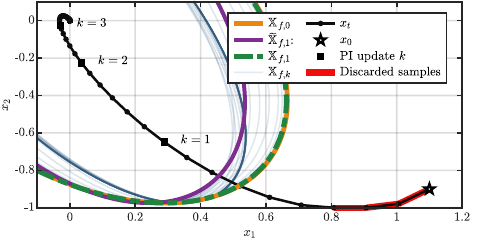}
    \caption{MPC policy iteration with damping and immediate update of the terminal pair. Here, $\alpha_1 = 0.0281$.}
    \label{fig:mpc_damp}
\end{figure}

\begin{figure}[t]
    \centering
    \includegraphics[width=\linewidth]{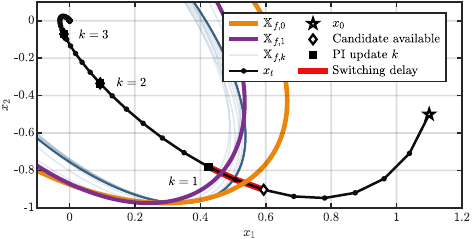}
    \caption{MPC policy iteration with reduced damping and delayed switching. Here, $\alpha_1=0.64$.}
    \label{fig:mpc_wait}
\end{figure}

Lastly, we demonstrate the closed-loop behavior of the proposed constrained MPC PI algorithm on the same system. The terminal cost matrix is initialized with \(P_0 = P_{\max}\). We choose $N=10$ and $M=2$, and use $H \geq 6$. We compare the two proposed strategies for safe updates of the terminal ingredients: damping and waiting, shown in Fig.~\ref{fig:mpc_damp} and Fig.~\ref{fig:mpc_wait} respectively.
The damping strategy enforces the terminal pair and immediate
switching conditions by reducing the step size.
Here, this results in a reduced step size of $\alpha_1 = 0.0281$ for the first update.
In contrast, the state-dependent waiting strategy enables a larger step size of $\alpha_1 = 0.64$, but delays the update until a feasible switching state is reached.
The damped update reaches $\alpha = 1$ at $k=9$ updates, while the waiting strategy reaches $\alpha = 1$ at $k=2$.
Subsequent updates are displayed in blue and converge to the 
invariant set associated with $P^\star$.

\section{Conclusions}
We proposed a policy iteration framework for receding-horizon control and model predictive control that updates the terminal cost using finite closed-loop data, with the terminal cost serving as both a policy parameter and a value function approximation.
In particular, we established a connection between the contraction properties of the Riccati map in the Thompson metric and approximate PI for RHC and MPC for the deterministic LQR problem. 
By constructing an invariant order interval, we have shown that the truncated PI step for RHC can be interpreted as a perturbation of a contractive Riccati map, yielding an ISS-type contraction bound with computable conditions for exact convergence. We have further developed a constructive extension to constrained MPC, in which the terminal cost and the associated terminal set are updated while preserving recursive feasibility and asymptotic closed-loop stability throughout the online updates.

While the LQR problem admits a well-known analytical solution, which enables the analysis of the effects of finite planning and rollout horizons on the contraction properties of the PI scheme, the proposed framework provides a new perspective on PI in the value-function space and its connection to predictive control. This work establishes a theoretical foundation for online terminal cost learning in RHC and MPC and opens up several avenues for future research, including extensions to nonlinear, robust, and stochastic settings.

\appendices

\section*{References}
\bibliographystyle{IEEEtran}
\bibliography{references}

\begin{IEEEbiography}[{\includegraphics[width=1in,clip]{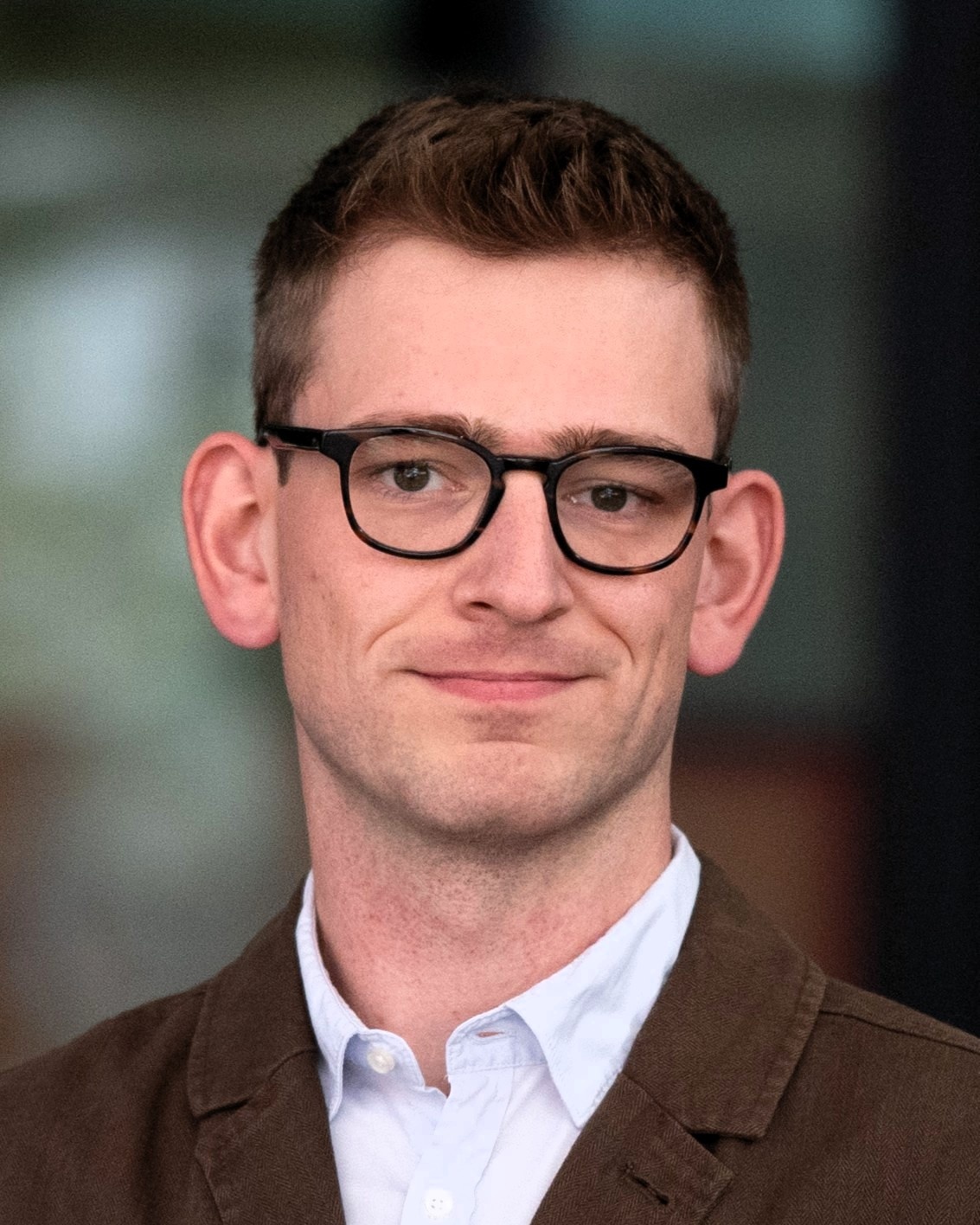}}]{Tren M.J.T. Baltussen} (Graduate Student Member, IEEE) received the M.Sc. degree in control theory (cum laude) from the Eindhoven University of Technology (TU/e), Eindhoven, The Netherlands, in 2024. He is currently pursuing the Ph.D. degree with the Control Systems Technology section, Department of Mechanical Engineering, TU/e. He held visiting researcher positions at The Ohio State University and at the University of California, Berkeley. His research interests include model predictive control, adaptive control, and learning-based control.
\end{IEEEbiography}

\begin{IEEEbiography}[{\includegraphics[width=1in,height=1.25in,clip,keepaspectratio]{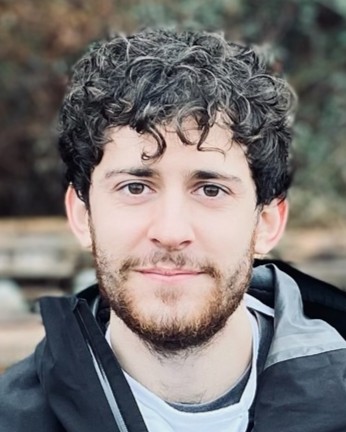}}]{Nathan P. Lawrence} (Member, IEEE) is a Senior Reinforcement Learning Engineer at HammerheadAI, Inc. He earned his Ph.D. in mathematics from the University of British Columbia in 2023 and was a Postdoctoral Scholar at the University of California, Berkeley before transitioning to industry in 2026. His research is centered around the safe design of reinforcement learning systems and sits broadly in optimization, model predictive control, and machine learning.
\end{IEEEbiography}

\begin{IEEEbiography}[{\includegraphics[width=1in,height=1.25in,clip,keepaspectratio]{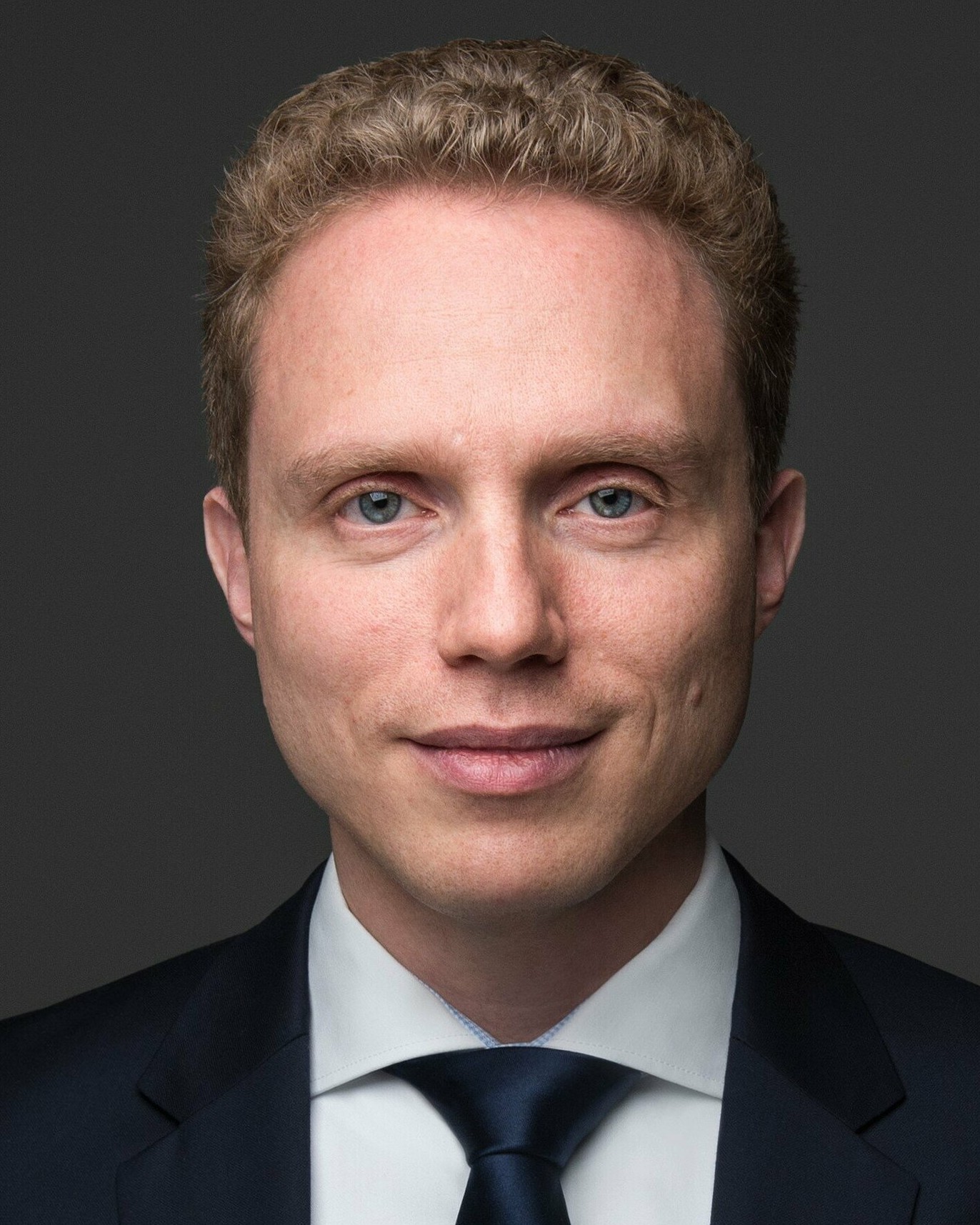}}]{Alexander Katriniok} (Senior Member, IEEE)
received the Ph.D. degree in mechanical engineering from RWTH Aachen University, Aachen, Germany, in 2013. He is currently an Assistant Professor with the Control Systems Technology section, Eindhoven University of Technology (TU/e), Eindhoven, The Netherlands. Dr. Katriniok is the Chair of the IEEE Control Systems Society (CSS) Technical Committee on Automotive Controls. He serves as an Associate Editor for IEEE Transactions on Control Systems Technology, the IEEE CSS Technology Conference Editorial Board, and the IEEE Robotics and Automation Letters. His research interests include learning-based control, provable safety and stability guarantees, and real-time numerical optimization, with application to (networked) autonomous systems.
\end{IEEEbiography}

\begin{IEEEbiography}[{\includegraphics[width=1in,height=1.25in,clip,keepaspectratio]{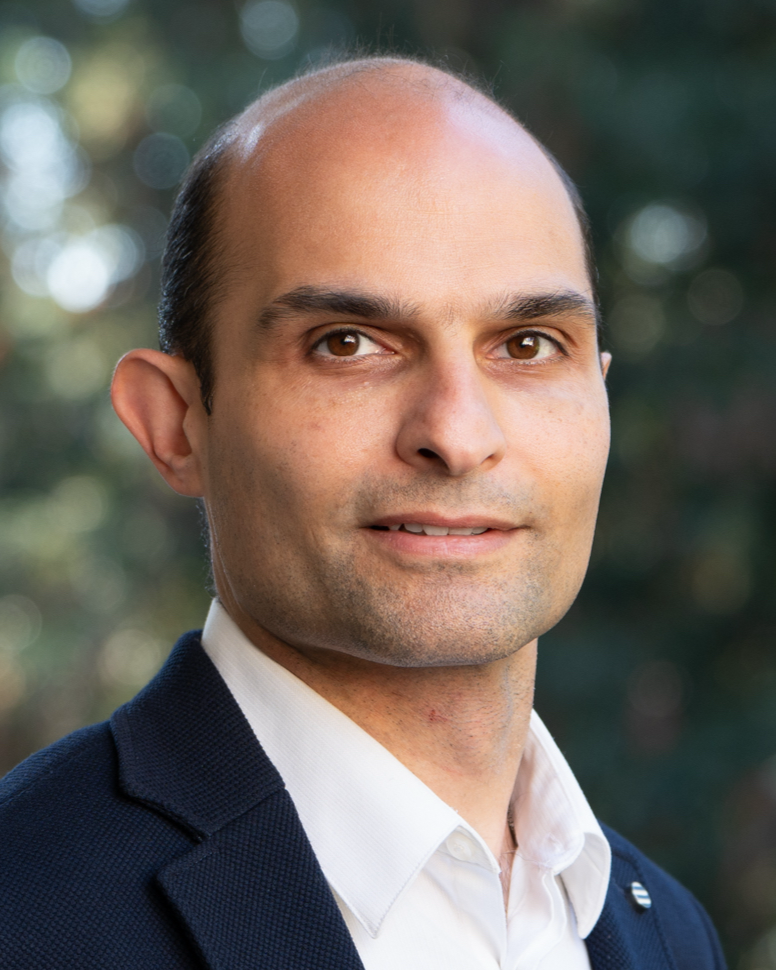}}]{Ali Mesbah} (Senior Member, IEEE) is Associate Professor of Chemical and Biomolecular Engineering at the University of California, Berkeley. He is a Senior Member of IEEE and AIChE. His research lies at the intersection of optimal control, machine learning, and applied mathematics, with applications to learning-based analysis, optimization, and predictive control of materials processing and biomedical systems.
\end{IEEEbiography}

\begin{IEEEbiography}[{\includegraphics[width=1in,height=1.25in,clip,keepaspectratio]{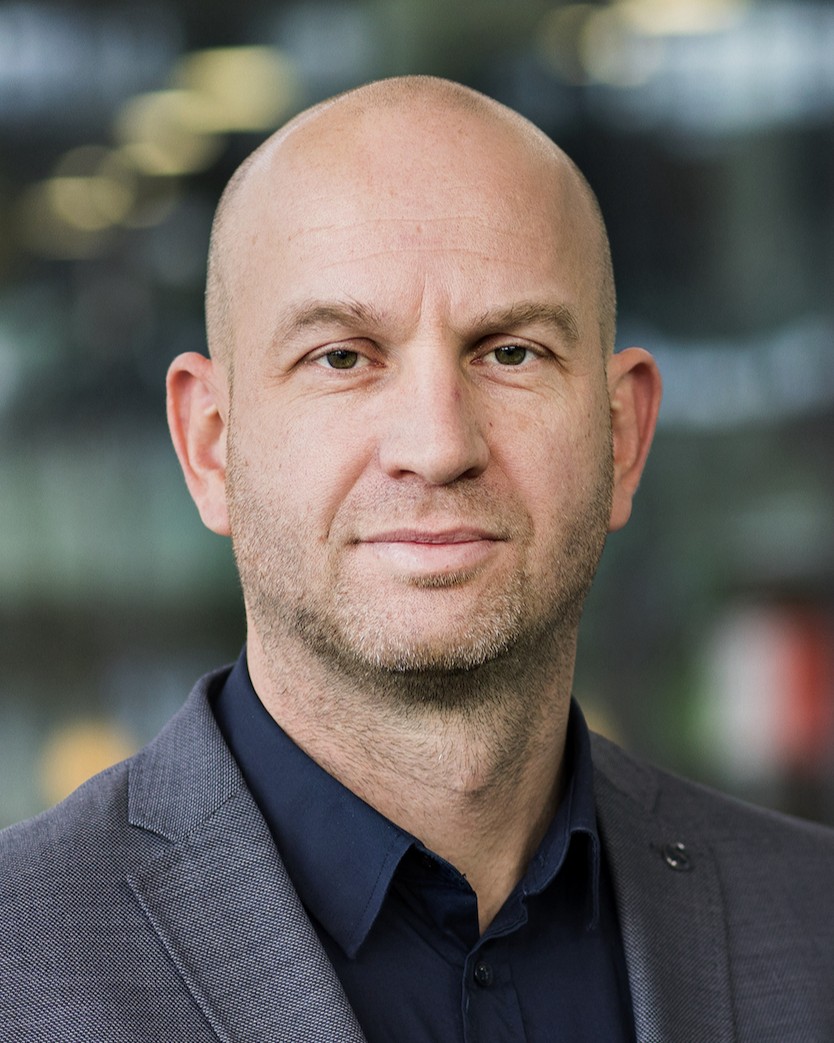}}]{W.P.M.H. (Maurice) Heemels} (Fellow, IEEE) received M.Sc. (mathematics) and Ph.D. (EE, control theory) degrees (summa cum laude) from the Eindhoven University of Technology (TU/e) in 1995 and 1999, respectively. From 2000 to 2004, he was with the Electrical Engineering Department, TU/e, as an assistant professor, and from 2004 to 2006 with TNO-Embedded Systems Institute as a Research Fellow. Since 2006, he has been with the Department of Mechanical Engineering, TU/e, where he is currently a Full Professor and Vice-Dean. He held visiting professor positions at ETH, Switzerland (2001), UCSB, USA (2008) and University of Lorraine, France (2020). He is a Fellow of IEEE, EUCA and IFAC, and was the chair of the IFAC Technical Committee on Networked Systems (2017-2023). He served/s on the editorial boards of Automatica,  Nonlinear Analysis: Hybrid Systems (NAHS), Annual Reviews in Control, and IEEE Transactions on Automatic Control, and is the Editor-in-Chief of NAHS as of 2023. He was a recipient of a personal VICI grant awarded by NWO (Dutch Research Council) and an ERC Advanced Grant. 
He was the recipient of the 2019 IEEE L-CSS Outstanding Paper Award, the Automatica Paper Prize 2020-2022 and an IFAC Outstanding Service Award in 2026. He was elected for the IEEE-CSS Board of Governors (2021-2023). His current research includes hybrid and cyber-physical systems, networked and event-triggered control systems and model predictive control and their applications.
\end{IEEEbiography}

\end{document}